\documentclass[11pt,english]{smfart}
\usepackage[applemac]{inputenc}
\def\today{janvier 2009} 
\usepackage{amsmath,amsfonts,amsthm,amssymb,amscd}
\usepackage[francais,english]{babel}
\newcommand{\be}{\begin{equation}}
\newcommand{\ee}{\end{equation}}
\newcommand{\bi}{\begin{itemize}}
\newcommand{\ei}{\end{itemize}}

\theoremstyle{plain} \newtheorem{theorem}{Theorem}[section]
\newtheorem{lemma}[theorem]{Lemma}
\newtheorem{proposition}[theorem]{Proposition}

\newtheorem{definition}[theorem]{Definition} \theoremstyle{remark}
\newtheorem{remark}[theorem]{Remark}

\newcommand{\R}{{\mathbb R}} \newcommand{\U}{{\mathcal U}}
 
\newcommand{\tHs}{{\tilde{H^s}}} 
 
\newcommand{\Z}{{\mathbb Z}}
\newcommand{\N}{{\mathbb N}}
\newcommand{\Nd}{{\mathbb N_d}}
\newcommand{\e}{{\varepsilon}}
\newcommand{\ba}{{\beta}}

\newcommand{\Zb}{{\bar \Z}} \newcommand{\Nb}{{\bar \N}}

\newcommand{\Tc}{{\mathcal T}}
\newcommand{\V}{{\mathcal V}}
\newcommand{\W}{{\mathcal W}}
\newcommand{\Ps}{{\mathcal P_{s}}}

\newcommand{\ls}{{\mathcal L_{s}}}
\newcommand{\Hs}{{\mathcal H^s}}
\newcommand{\Hsd}{{\mathcal H^s_d}}
\newcommand{\Qs}{{\mathcal Q_s}}

\newcommand{\C}{\mathbb{C}}
\newcommand{\T}{\mathcal{T}}
\newcommand{\Tt}{\mathbb{T}}

\newcommand{\va}[1]{|#1|}
  
  \newcommand{\vaj}{|j|}
  \newcommand{\val}{|l|}
  
  \newcommand{\vat}{|t|}
  \newcommand{\Va}[1]{\left|#1\right|}

\def\norma#1{\left\| #1\right\|}

\def\sleq{\leq\kern-6pt \cdot\null\hskip4pt}

\numberwithin{equation}{section}

\begin{document}

\author{ Beno\^it Gr\'ebert, Rafik Imekraz, Eric Paturel} \title[Normal forms for semilinear quantum harmonic oscillators ]{Normal forms for semilinear quantum harmonic oscillators}
\alttitle{Formes normales de Birkhoff pour l'oscillateur harmonique quantique non lin\'eaire}
\date{\today}

\begin{abstract}
We consider the semilinear harmonic oscillator
$$i\psi_t=(-\Delta +\va{x}^{2} +M)\psi +\partial_2 g(\psi,\bar \psi), \quad x\in \R^d,\  t\in \R$$
where $M$ is a Hermite multiplier and $g$ a smooth function globally of order 3 at least.\\
We prove that such a Hamiltonian equation admits, in a neighborhood of the origin, a Birkhoff normal form at any order and that, under generic conditions on $M$ related to the non resonance of the linear part, this normal form is integrable when $d=1$ and gives rise to  simple (in particular bounded) dynamics when $d\geq 2$.\\
As a consequence we prove the almost global existence  for solutions of the above equation with small Cauchy data. Furthermore we control the  high Sobolev norms of these solutions.
    \end{abstract}

\begin{altabstract} Dans cet article nous consid\'erons l'oscillateur harmonique semi-lin\'eaire:
$$i\psi_t=(-\Delta +x^{2} +M)\psi +\partial_2 g(\psi,\bar \psi), \quad x\in \R^d,\  t\in \R$$
o\`u $M$ est un multiplicateur de Hermite  et $g$ est une fonction r\'eguli\`ere globalement d'ordre au moins trois.\\
Nous montrons qu'une telle \'equation admet, au voisinage de z\'ero, une forme normale de Birkhoff \`a n'importe quel ordre et que, sous des hypoth\`eses g\'en\'eriques sur $M$ li\'ees \`a la non r\'esonance de la partie lin\'eaire, cette forme normale est compl\`etement int\'egrable si $d=1$ et donne lieu \`a une dynamique simple (et en particulier born\'ee) pour $d\geq 2$.\\
Ce r\'esultat nous permet de d\'emontrer l'existence presque globale et de contr\^oler les normes de Sobolev d'indice grand des solutions de l'\'equation non lin\'eaire ci-dessus avec donn\'ee initiale petite.
  
\end{altabstract}    
\altkeywords{Forme Normale de Birkhoff, Osillateur harmonique non lin\'eaire, 
EDP Hamiltoni\`ennes, stabilit\'e pour des temps longs, \'equation de Gross-Pitaevski}\frontmatter
\maketitle
\newpage
\mainmatter

\section{Introduction, statement of the results}
The aim of this paper is to prove a Birkhoff normal form theorem  for the semilinear harmonic oscillator equation
\begin{align}\Big\{ \begin{split}
\label{1}
 i\psi_t&=(-\Delta +\va{x}^{2} +M)\psi +\partial_2 g(\psi,\bar \psi)\\ 
\psi\arrowvert_{t=0}&=\psi_0 
\end{split}\end{align}
on the whole space $\R^d$ ($d\geq 1$) and  to discuss its dynamical consequences. Here  $g$ is a smooth function, globally of order $p\geq 3$ at 0, and $\partial_2 g$ denotes the partial derivative of $g$ with respect to the second variable. The linear operator $M$ is a Hermite multiplier. To define it precisely (at least in the case $d=1$, see Section \ref{nlsdd} for the multidimensional case), let us introduce
the quantum harmonic oscillator on $\R^d$, denoted by $T=-\Delta +\va{x}^2$. When $d=1$, $T$  is diagonal in the Hermite basis $(\phi_j)_{j\in\Nb}$:
\begin{align*}T\phi_j &=(2j-1)\phi_j, \quad j\in \Nb\\
\phi_{n+1} &=\frac{H_n(x)}{\sqrt{2^n n!}}e^{-x^2/2} , \quad n\in \N
\end{align*}
where $H_n(x)$ is the $n^{th}$ Hermite polynomial relative to the weight $e^{-x^2}$ :
$$
\int_\R e^{-x^2} H_m(x)H_n(x)dx=2^n n! \sqrt\pi \delta_{nm}.$$
In this basis (and for $d=1$), a Hermite multiplier is an operator given by
$$
M\phi_j=m_j \phi_j\,,
$$
where $(m_j)_{j\in \Nb}$ is a bounded sequence of real numbers, that will be chosen in the following classes : for any $k\geq 1$, we define the class
$$
\W_k=\{   (m_j)_{j\in \Nb}\mid \mbox{ for each }j,\ m_j=\frac{\tilde m_j}{j^k}\mbox{ with } \tilde m_j\in [-1/2,1/2]\}
$$
that we endow with the product probability measure. In this context the linear frequencies, i.e. the eigenvalues of $T+M= -d^2/dx^2 +x^2+M$  are given by
\begin{equation} \label{ome}\omega_j= 2j-1+m_j=2j-1 +\frac{\tilde m_j}{j^k}, \quad j\in \Nb.\end{equation}

Let 
\begin{align}\begin{split} \label{H}
\tilde{H^s}=\{&f\in H^s(\R^d,\C) | x\mapsto {x}^{\alpha}\partial^\ba f \in L^2(\R^d)\\
&\mbox{ for any } \alpha,\  \beta\in \N^d \mbox{ satisfying } 0\leq \va{\alpha}+\va{\ba}\leq s \}
\end{split}\end{align}
where $H^s(\R^d,\C)$ is the standard Sobolev space on $\R^d$. We note that, for any $s\geq 0$,
the domain of $T^{s/2}$ is $\tHs$ (see for instance \cite{Helf84} Proposition 1.6.6) and that for $s>d/2$, $\tHs$ is an algebra.\\
If $\psi_0\in \tilde{H^s}$ is small, say of norm $\epsilon$, local existence theory implies that \eqref{1} admits a unique solution in $\tilde{H^s}$ defined on an interval of length $c\epsilon^{-p+2}$. Our goal is to prove that for $M$ outside an exceptional subset, given any integer $r\geq 1$ and provided that $s$ is large enough and $\epsilon$ is small enough,  the solution extends over an interval of length $c\epsilon^{-r}$. Furthermore we  control the norm of the solution in $\tilde{H^s}$-norm ($d\geq 1)$ and localize the solution in the neighborhood of a torus (only in the case $d=1$, cf. Theorem \ref{thm:dyn} and Theorem \ref{thm:dynd}).\\
Precisely  we have
\begin{theorem} \label{thm1} 
Let $r,k\in \N$ be  arbitrary integers. There exists a set $F_{k}\subset \W_{k}$
whose measure equals $1$ such that if $m=(m_j)_{j\in \Nb}\in F_{k}$ and if  $g$ is a $ C^\infty$ function on a neighborhood of the origin in $\C^2$, satisfying $g(z,\bar z)\in \R$ and
vanishing at least at order 3 at the origin, there is $s_{0}\in\N$ such that  
for any $s\geq s_{0}$, f there are $\epsilon_{0}>0$, 
$c>0$, such that for any 
$\epsilon \in (0,\epsilon_{0})$,  for any  $\psi_0$ in $\tHs$ with $||\psi_0||_s \leq \epsilon$ ,  the Cauchy problem \eqref{1} with initial datum $\psi_0$ 
has a unique solution 
$$ \psi\in 
C^1((-T_{\epsilon},T_{\epsilon}),\tilde H^{s})$$ with
$T_{\epsilon}\geq c\epsilon^{-r}$. Moreover,
for any $t\in (-T_{\epsilon},T_{\epsilon})$, one has
\be
\label{estimHs} \Vert{\psi(t,\cdot )}\Vert_{\tilde H^{s}}\leq 2\epsilon \ .  
\ee
\end{theorem}
For the nonlinearity $g(\psi,\bar \psi)=\lambda \frac 2 {p+1} \va{\psi}^{p+1}$ with $p\geq 1$ and without Hermite multiplier ($M=0$), we recover the Gross-Pitaevskii equation
\be\label{GP}
i\psi_t=(-\Delta +\va{x}^{2} )\psi +\lambda \va{\psi}^{p-1}\psi, \quad t\in \R, \quad x\in \R^d.
\ee
In this case, the global existence  in the energy space  $\tilde H^{1}$ has been proved  for\footnote{we use the convention  $\frac{d+2}{(d-2)^+}=+\infty$ for $d=1,2$, and $(d-2)^+ =d-2$ for $d\geq 3$} $1\leq p<\frac{d+2}{(d-2)^+}$ without smallness assumption on the Cauchy data in the defocusing case ($\lambda <0$) and for small Cauchy data in the focusing case ($\lambda>0$) (see \cite{Car02} and also
\cite{Zha05}). But nothing is known for nonlinearities of higher order, neither about conservation of the $\tilde H^{s}$-norm for $s>1$. Our result states that, avoiding resonances by adding a generic linear term $M\psi$ (but $M=0$ is not allowed), we recover almost global existence for solutions of Gross-Pitaevskii equation with a nonlinearity of arbitrary high order and small Cauchy data in $\tilde H^{s}$ for $s$ large enough. In some sense, this shows that the instability for Gross-Pitaevskii that could appear in that regime are necessarily produced by resonances. More precisely, we can compare with the semi-classical cubic Gross-Pitaevskii in $\R^3$ which appears in the study of Bose-Einstein condensates (for a physical presentation see \cite{PS03})
\be\label{GPSC}
ihu_t=-h^2\Delta u+\va{x}^{2} u +h^2\va{u}^{2}u, \quad t\in \R, \quad x\in \R^3
\ee
where $h$ is a small parameter.\\
The scaling relation between $\psi$ solution of \eqref{GP} and $u$ solution of \eqref{GPSC} is given by
\be
u(t,x)=\frac{1}{\sqrt{h}}\psi(t,\frac{x}{\sqrt h}).
\ee
We note that for  multi indices $\alpha, \beta \in \N^3$, with $y = \frac{x}{\sqrt{h}}$, 
$$
\norma{y^\alpha \partial^\beta \psi}^2_{L^2(\R^3)}=h^{\va{\beta}-\va{\alpha}-1/2} \norma{x^\alpha \partial^\beta u}^2_{L^2(\R^3)}.$$
Thus the smallness of $\psi_0$ in $\tilde H^s$ imposed in Theorem \ref{thm1}, i.e. $ \Vert{\psi_0}\Vert_{\tilde H^{s}}\leq C\epsilon$, reads 
$$
\sum_{\va{\beta}+\va{\alpha}\leq s}h^{\va{\beta}-\va{\alpha}-1/2} \norma{x^\alpha \partial^\beta u_0}^2_{L^2(\R^3)}\leq C\epsilon^2.$$
Taking $\epsilon= h^{1/6}$ with $h$ small enough, this allows the derivatives of order greater than 1 to have large $L^2$-norm when $h$ is small:
$$\norma{ \partial^\beta u_0}^2_{L^2(\R^3)}=O(h^{-\va{\beta}+5/6})$$
i.e. the initial data has to be small in $L^2$ but may have large oscillations. Then, Theorem \ref{thm1} states that, avoiding the resonances by adding a generic linear term  (which, in the preceding scaling, stays of order $h$), the same estimates remain true for the solution $u(t,.)$ with $\vat=O(h^{-r/6})$, $r$ being chosen arbitrarily from the principle.  Notice that the role of the linear operator $M$  is  to remove the resonances between the free modes (see \eqref{ome}). The fully resonant case $M=0$ is beyond the scope of the paper.

\medskip

To prove Theorem \ref{thm1} we use  the Birkhoff normal form theory.
This technique  has been developped  by Bourgain \cite{Bo96}, Bambusi \cite{Bam03}, Bambusi-Grébert \cite{BG04} for semilinear PDEs (typically semilinear Schr\"odinger equation or semilinear wave equation) on the one dimensional torus and by Bambusi-Delort-Grébert-Szeftel \cite{BDGS} for the semilinear Klein-Gordon equation on the sphere $S^d$ (or a Zoll manifold). These cases were concerned with compact domains. In our work the domain is $\R^d$, the potential $x^2$ guarantees that the spectrum remains pure point, but the free modes of the harmonic oscillator are not so well localized.\\
For general reference on Hamiltonian PDEs and their perturbations, see the recent monographies \cite{Cr00, K2, Bo05, KaP}. We also note that in \cite{K1}, a KAM-like theorem is proved for \eqref{1} in one dimension and with special nonlinearities .

Let us describe roughly the general method. Consider a Hamiltonian system whose Hamiltonian function decomposes in a quadratic part, $H_0$ (associated to the linear part of the equation), and a perturbative nonlinear part $P$ (at least cubic): $H=H_0+P$. We assume that $H_0$ is diagonal in a Hilbert basis $(\phi_j)_{j\geq 1}$ of the phase space $\mathcal P$ : $H_0=\sum_j
\omega_j \xi_j \eta_j$ for $(\xi,\eta)\in \mathcal P$  and $\omega=(\omega_j)_{j\geq 1}$ is the vector of free frequencies (the eigenvalues of the linear part). In the harmonic oscillator case, the Hilbert basis is given by the Hermite functions and $\mathcal P =\ell^2\times \ell^2$. The heuristic idea could be resumed as follows: if the free modes do not interact linearly  (i.e. if  $\omega$ is non resonant), and if they do not interact too much via the nonlinear term, then the system will remain close to an integrable one, up to a nonlinear term of very high order, and thus the solutions will exist and stay under control during a very long time. More precisely, by a Birkhoff normal form approach we prove  (cf. Theorem \ref{thm:birk} which is our main theorem) that  $H\sim H'_0+P'$ where $H'_0$ is no more quadratic but remains integrable (in the case $d=1$) and $P'$ is at least of order $r$, where $r$ can be chosen arbitrarily large as soon as we work in a sufficiently small neighborhood of the origin.\\
 To guarantee the second condition, i.e. that the free modes do not interact too much via the nonlinear term, we have to control the integral of the product of three or more modes:
\be\label{aj}a_j=\int_D \phi_{j_1}(x)\cdots \phi_{j_k}(x) dx\ee
where $D$ is the space domain ($\R^d$ in our case) and $j$ is a multi-index in $\N^k$, $k$ being smaller than the fixed order $r$ and larger than $3$. It turns out that, in our case, this control cannot be as good as in the  cases of compact domains studied previously.\\
 Let us consider ordered multi-indices $j$, i.e. such that $j_1\geq j_2\geq \cdots \geq j_k$.
In \cite{BDGS, Gre07,Bam07}   the following control was used: there exists $\nu>0$ and for any $N\geq 1$ there exists $C_N>0$ such that for all ordered $j$
\be\label{aj1}
\va{a_j}\leq C_N j_3^\nu \left( \frac{j_3}{j_3+j_1-j_2}\right)^N.
\ee
In the case of the harmonic oscillator, this estimate is false (cf \cite{Wa08} where an equivalent is computed for four modes) and we are only able to prove the following: there exists $\nu>0$ and
for any $N\geq 1$ there exists $C_N>0$ such that for all ordered $j$
\be\label{aj2}
\va{a_j}\leq C_N \frac{j_3^\nu}{j_1^{1/24}} \left( \frac{\sqrt{j_2j_3}}{\sqrt{j_2j_3}+j_1-j_2}\right)^N.
\ee
The difference could seem minimal but it is technically important: \\
$
\sum_{j_1}\left( \frac{j_3}{j_3+j_1-j_2}\right)^\mu\sim C j_3$ for an uniform constant $C$ providing $\mu>1$ and similarly $\sum_{j_1}\left( \frac{\sqrt{j_2j_3}}{\sqrt{j_2j_3}+j_1-j_2}\right)^\mu\sim C\sqrt{j_2j_3}$ for $\mu>1$.  In the first case, the extra term $j_3$ can be absorbed by changing the value of $\nu$ in \eqref{aj1} ($\nu'=\nu +1$). This is not possible in the second case. In some sense the perturbative nonlinearity is no longer short range (cf. \cite{Wa08}).\\
Actually in the case studied in \cite{Bo96, Bam03, BG04}, the linear modes (i.e. the eigenfunctions of the linear part) are localized around the exponentials $e^{ikx}$, i.e. the eigenfunctions of the Laplacian on the torus. In particular the product of eigenfunctions is close to an other eigenfunction which makes the control of \eqref{aj}  simpler. In the harmonic oscillator case, the eigenfunctions are not localized and the product of eigenfunctions has more complicated properties. Notice that, in the case of the semilinear Klein-Gordon equation on the sphere, the control of \eqref{aj} is more complicated to obtain,  but an estimate of type \eqref{aj1} is proved in \cite{DS1} for the Klein-Gordon equation on Zoll manifolds.\\
From the point of view of a normal form, the substitution of \eqref{aj1} by \eqref{aj2} has the following consequence:\\
Consider a formal polynomial
$$Q(\xi, \eta) \equiv Q(z)=\sum_{l=0}^k\sum_{j\in 
\N^l}a_{j}z_{j_{1}}\ldots z_{j_{l}}$$ 
with coefficients $a_j$ satisfying \eqref{aj1}. In \cite{Gre07} or \cite{Bam07}, it is proved that its Hamiltonian vector field $X_Q$ is then regular from\footnote{here $l^2_s=\{(z_l)\mid \sum l^{2s}\va{z_l}^2<\infty\}$ and corresponds to functions $\psi=\sum z_l\phi_l$ in $\tilde{H^{2s}}$.} $\Ps=\ell^2_s\times \ell^2_s$
 to $\Ps$ for all $s$ large enough (depending on $\nu$). In our present case, i.e. if $a_j$ only satisfy \eqref{aj2}, which defines the class $\Tc^\nu$, then we prove that $X_Q$ is regular from $\Ps$ to $\mathcal P_{s'}$   for all $s'<s-1/2+1/24$ and $s$ large enough. This "loss of regularity" would  of course complicate an iterative procedure, but it is bypassed in the following way : the nonlinearity $P$ is regular in the sense that $X_P$ maps  $\Ps$ to $\mathcal P_{s}$ continuously for $s$ large enough (essentially because the space $\tilde H^s$ is an algebra for $s>d/2$). On the other hand, we build at each step a canonical transform which preserves the regularity. Indeed, at each iteration, we compute the canonical transformation as the time 1 flow of a Hamiltonian $\chi$, and the solution of the so called homological equation gives rise to an extra term in \eqref{aj2} for the coefficient of the polynomial $\chi$:
\be\label{aj3}
\va{a_j}\leq C_N \frac{j_3^\nu}{j_1^{1/24}(1+j_1-j_2)} \left( \frac{\sqrt{j_2j_3}}{\sqrt{j_2j_3}+j_1-j_2}\right)^N.
\ee 
Using such an estimate on the coefficients (in the class\footnote{Actually in section \ref{poly}, instead of $\Tc^{\nu}$ and $\Tc^{\nu,+}$, we consider  more general classes $\Tc^{\nu,\beta}$ and $\Tc^{\nu,\beta,+}$ where the parameter $\beta$ plays the role of the exponent $1/24$ in \eqref{aj2} and \eqref{aj3}} denoted $\Tc^{\nu,+}$ in Section \ref{poly}), we prove in Proposition \ref{tame} that $X_\chi$
is regular from $\Ps$ to $\Ps$   for all $s$ large enough. Furthermore, we prove in Proposition \ref{bra} that the Poisson bracket of a polynomial in $\Tc^\nu$ with a polynomial in $\Tc^{\nu,+}$ is in $\Tc^{\nu'}$ for some $\nu'$ larger than $\nu$. So an iterative procedure is possible in $\Ps$. \\
This smoothing effect of the homological equations was already used by S. Kuksin in \cite{Kuk87} (see also \cite{K1, Pos96}). Notice that this is, in some sense, similar to the local smoothing property for Schr\"odinger equations with potentials superquadratic at infinity  studied in \cite{YZ04}.

\medskip

Our article is organized as follows: in Section \ref{birk} we state and prove a specific Birkhoff normal form theorem adapted to the loss of regularity that we explained above. In Section \ref{dyn}, we apply this theorem to the $1-d$ semilinear harmonic oscillator equation (Subsection \ref{nls1d}) and  we generalize it to cover the multidimensional case (Subsection \ref{nlsdd}).

 {\it Acknowledgements: it is a great pleasure  to
thank Dario Bambusi and Didier Robert for many helpful discussions. We thank both referees for useful suggestions.}

\section{The Birkhoff normal form}\label{birk}
\subsection{The abstract model}\label{model}
To begin with, we give an abstract model of infinite dimensional
Hamiltonian system. In Section \ref{dyn} we will verify that the  
nonlinear harmonic oscillator can be described in this abstract framework. Throughout the paper, we denote $\Nb = \N  
\setminus \{0\}$ and $ \Zb = \Z \setminus \{0\}$.
We work in the phase space $\Ps\equiv\Ps(\C):=\ell^2_{s}(\C)\times  
\ell^2_{s}(\C)$
where, for $s\in \R_+$,
$\ell^2_{s}(\C):=\{(a_{j})_{j\geq 1}\in \C^\Nb \mid \sum_{j\geq
1}j^{2s}\va{a_{j}}^2 < +\infty  \}$  is a Hilbert space for the
standard norm: $\norma{a}_{s}^2=\sum_{j\geq
1}\vaj^{2s}\va{a_{j}}^2$.  We denote  
$\Ps(\R):=\{(\xi,\bar\xi)\in\Ps(\C)\}$ the "real" part of $\Ps(\C)$. We shall denote a general point of $\Ps$ by $z=(\xi, \eta)$ with $z=(z_j)_{j\in \Zb}$, $\xi=(\xi_j)_{j\in \Nb}$, $\eta=(\eta_j)_{j\in \Nb}$ and the correspondence: $z_j=\xi_j, \ z_{-j}=\eta_j$ for all $j\in \Nb$. Finally, for a Hamiltonian function $H$, the Hamiltonian vector field $X_H$ is defined by $$X_H(z) = \left( \left(- \frac{\partial H}{\partial \xi_k}\right)_{k \in \Nb}, \left(\frac{\partial H}{\partial \eta_k}\right)_{k \in \Nb} \right).$$

\begin{definition} \label{defHs}Let $s\geq 0$, we denote by
$\Hs$ the space of Hamiltonian functions $H$ defined on
a neighborhood $\U$ of the origin in $\Ps\equiv\Ps(\C)$,
satisfying $H(\xi,\bar \xi)\in \R$ (we say that $H$ is real) and
$$H\in C^\infty (\U, \C)\quad \mbox{and} \quad X_{H}\in C^\infty (\U,  
\Ps),$$ as well as every homogeneous polynomial $H_k$ appearing in the Taylor expansion of $H$ at $0$ :
$$H_k\in C^\infty (\U, \C)\quad \mbox{and} \quad X_{H_k}\in C^\infty (\U,  
\Ps).$$

\end{definition}
\begin{remark}
This property, for Hamiltonians contributing to the nonlinearity, will in particular force them to be  {\em semilinear perturbations} of the harmonic oscillator.
\end{remark}
In particular the Hamiltonian vector fields of functions $F,\ G$ in
$\Hs$ are in $\ell^2_{s}(\C)\times \ell^2_{s}(\C)$ and we can define  
their Poisson
bracket by
$$
\{ F,G \}=i\sum_{j\geq 1} \ \frac{\partial F}{\partial \xi_{j}}
\frac{\partial G}{\partial \eta_{j}}-
\frac{\partial F}{\partial \eta_{j}}\frac{\partial G}{\partial \xi_{j}}
\ .$$
Notice that since for $P\in\Hs$, the vector field $X_P$ is a $C^\infty$ function from a neighborhood of $\Ps$ to $\Ps$ we have
\begin{lemma}\label{first}
Let  $P \in \Hs$ such that $P$ vanishes up to order $r+1$ at the origin, 
that is :
$$ \forall k \leq r+1, \forall j \in \Zb^k, \; \frac{\partial^k P}{\partial z_{j_1} \dots \partial z_{j_k}}(0) = 0$$
Then there exists $\varepsilon_0 >0$ and $C>0$ such that, for  $z\in \Ps $ satisfying 
$||z||_s \leq \varepsilon_0$, we have
$$ ||X_P(z)||_s \leq C ||z||_s^r. $$
\end{lemma}

Our model of integrable system is the harmonic oscillator
$$H_{0}=\sum_{j\geq 1} \omega_{j}\xi_{j}\eta_{j}$$
where $\omega =(\omega_{j})_{j\geq 1}\in \R^\Nb$ is the frequency
vector. We will assume that these frequencies grow at most polynomially,
i.e. that there exist $C>0$ and  $\bar{d}\geq 0$ such that for any  
$j\in \Nb$,
  \be \label{om}|\omega_{j}|\leq C \vaj^{\bar{d}},\ee
  in such a
way that  $H_{0}$ be well defined on $\Ps$ for $s$ large enough.\\
The perturbation term is a real  function, $P\in
\Hs$, having a zero of order at least $3$ at the
origin. Our Hamiltonian function is then given by
$$H=H_{0}+P$$
and Hamilton's canonical equations read
\be \label{HeqC}
\left\{\begin{array}{ccc}
\dot \xi_{j} &=&-i\omega_{j} \xi_{j} -i\frac{\partial P}{\partial  
\eta_{j}},\
j\geq 1\\
\dot \eta_{j} &=& i\omega_{j}\eta_{j}+i\frac{\partial P}{\partial  
\xi_{j}},\
j\geq 1.
\end{array}
\right.
\ee
Our theorem will require essentially two hypotheses: one on the
perturbation $P$ (see Definition \ref{po}) and one  
on the frequency vector
$\omega$ that we describe now.

For $j\in \Zb^k$ with $k\geq 3$, we define
$\mu(j)$ as the third largest integer among
$\va{j_{1}},\ldots,\va{j_{k}}$. Then we set
$S(j):=\va{j_{i_{1}}}-\va{j_{i_{2}}}$ where
$\va{j_{i_{1}}}$ and $\va{j_{i_{2}}}$ are respectively the largest  
integer and the second
largest integer among
$\va{j_{1}},\ldots,\va{j_{k}}$. In particular, if the multi-index $j$  
is ordered i.e.
if $
\va{j_{1}}\geq\ldots\geq\va{j_{k}} $ then
$$\mu(j):=\va{j_{3}} \mbox{ and }
S(j)=\va{j_{1}}-\va{j_{2}}.$$

In \cite{Bam03, BG04, Gre07,Bam07} the non resonance condition on  
$\omega$ reads
\begin{definition}\label{NR} A frequency vector $\omega \in \R^\Nb$
is \textbf{ non resonant} if for
     any $r\in\Nb$, there are $ \gamma >0$ and
     $\delta >0$ such that for any $j\in \Nb^{r}$ and any $1\leq i\leq
     r$, one has
     \begin{equation}
     \label{A.1}
     \left|\omega_{j_1}+\cdots+\omega_{j_{i}}-\omega_{j_{i+1}}-\cdots
     -\omega_{j_{r}} \right|\geq \frac{\gamma}{\mu (j)^{\delta}}
     \end{equation}
except in the case $\{j_{1},\ldots,j_{i}\}=\{j_{i+1},\ldots,j_{r}\}$.
\end{definition}
In the harmonic oscillator case\footnote{ The following holds, more generally, if the frequency vector is non resonant as in Definition \ref{NR}  
and satisfies the asymptotic: $\omega_l\sim l^n$ with $n\geq 1$. }, we  
are  able to work with a slightly refined non resonance condition
\begin{definition}\label{SNR} A frequency vector $\omega \in \R^\Nb$
is \textbf{strongly non resonant} if for
     any $r\in\Nb$, there are $ \gamma >0$ and
     $\delta >0$ such that for any $j\in \Nb^{r}$ and any $1\leq i\leq
     r$, one has
     \begin{equation}
     \label{A.2}
     \left|\omega_{j_1}+\cdots+\omega_{j_{i}}-\omega_{j_{i+1}}-\cdots
     -\omega_{j_{r}} \right|\geq \gamma \frac{1+S(j)}{\mu(j)^{\delta}}
     \end{equation}
except if $\{j_{1},\ldots,j_{i}\}=\{j_{i+1},\ldots,j_{r}\}$.
\end{definition}
This improvement of the non resonance condition is similar to the modification to the standard  second Melnikov condition introduced first by S. Kuksin in \cite{Kuk87} (see also \cite{K1} and \cite{Pos96}).

\subsection{Polynomial structure}\label{poly}
For $j\in \Zb^k$ with $k\geq 3$, we have already defined $\mu(j)$ and  
$S(j)$, we now introduce
$$B(j)={ \va{j_{i_{2}}\ j_{i_{3}}}}^{1/2}, \quad C(j)=\va{j_{i_{1}}}$$ where
$\va{j_{i_{1}}}$, $\va{j_{i_{2}}}$ and $\va{j_{i_{3}}}$ are respectively the first,  the second  
 and the third largest integer among
$\va{j_{1}},\ldots,\va{j_{k}}$. We also define
\be\label{A}
A(j)=\frac{B(j)}{B(j)+S(j)}.
\ee
In particular, if the multi-index $j$ is ordered i.e.
if $
\va{j_{1}}\geq\ldots\geq\va{j_{k}} $ then
$$A(j)=\frac{\va{j_2 j_3}^{1/2}}{\va{j_2 j_3}^{1/2}+\va{j_1}-\va{j_2}}$$
and 
$$C(j)=\va{j_1}.$$

\begin{definition}\label{po} Let $k\geq 3$, $\beta \in (0,+\infty)$ and  
$\nu \in [0,+\infty)$ and
let \be \label{defQ}Q(\xi, \eta) \equiv Q(z)=\sum_{j\in
\Zb^k}a_{j}z_{j_{1}}\ldots z_{j_{l}}\ee be a formal homogeneous polynomial of degree  
$k$
on $\Ps(\C)$. $Q$ is in the class
$\T_{k}^{\nu,\beta}$ if for any $N\geq 1$ there exists a constant  
$c_N>0$ such that for all
$j\in \Zb^k$
\be\label{a}\va{a_{j}}\leq c_N \frac{\mu(j)^{\nu}}{C(j)^\beta}  
A(j)^N.\ee
\end{definition}
We will also use
\begin{definition}\label{po+} Let $k\geq 3$, $\beta \in [0,+\infty)$  
and $\nu \in [0,+\infty)$ and
let $$Q(\xi, \eta) \equiv Q(z)=\sum_{j\in
\Zb^k}a_{j}z_{j_{1}}\ldots z_{j_{l}}$$ be a formal homogeneous polynomial of degree  
$k$
on $\Ps(\C)$. $Q$ is in the class
$\T_{k}^{\nu,\beta,+}$ if for any $N\geq 1$ there exists a constant  
$c_N>0$ such that for all
$j\in \Zb^k$
\be\label{a+}\va{a_{j}}\leq c_N \frac{\mu(j)^{\nu}}{C(j)^\beta(1+S(j))}  
A(j)^N.\ee
\end{definition}
The best constants $c_N$ in \eqref{a} define a family of semi-norms for  
which
$\T_{k}^{\nu,\beta}$ is a Fr\'echet space.
\begin{remark}
Notice that the formula \eqref{defQ} does not give a unique representation of polynomials on $\Ps$. However, since the estimates \eqref{a} and \eqref{a+}  are symmetric with respect to the order of the indexes $j_1,\cdots,j_k$,  this non uniqueness does not affect Definitions \ref{po} and \ref{po+}.
\end{remark}

\begin{remark}
In the estimate \eqref{a}, the numerator allows an increasing behaviour
with respect to $\mu(j)$ that will be useful to control the small
divisors. The denominator imposes a slightly decreasing behaviour with respect  
to the largest index $C(j)$ and a highly decreasing behaviour for
monomials having their two modes of
largest indexes that are not of the
same order. This control is slightly better in $\T_{k}^{\nu,\beta ,+}$.  
\end{remark}

\begin{remark}
We will see in Proposition \ref{tame} that, if $\beta>1/2$ then $\T_{k}^{\nu,\beta}\subset  
\Hs$
for $s\geq \nu +1$. Unfortunately $\beta$ is not that large in the harmonic oscillator case, where the best we obtain is $\beta = 1/24$. Thus  $P\in\T_{k}^{\nu,\beta}$ does not imply $
P\in \Hs$. Nevertheless, as we will see in Proposition \ref{tame}, a polynomial in $\T_{k}^{\nu,\beta}$ is well defined and continuous on a
neighborhood of the origin in $\Ps(\C)$ for $s$ large enough. As a comparison, in \cite{Gre07, Bam07}, our estimate \eqref{a} is replaced with
\be\label{a++}\va{a_{j}}\leq C_N  
\frac{\mu(j)^{N+\nu}}{(\mu(j)+S(j))^N}.\ee
which is actually better than (\ref{a}), since it implies the $\Hs$ regularity. This type of control on the coefficients $a_j$ was first introduced in \cite{DS1} in the context of multilinear forms.
\end{remark}


\begin{definition}\label{T} Let $\nu\geq 0$ and $\beta\geq 0$. A  
function $P$ is in the class
$\T^{\nu,\beta}$ if
\begin{itemize}
     \item
there exists $s_{0}\geq 0$ such that, for any $s\geq s_{0}$ there exists $\U_s$, a neighborhood of the origin in $\Ps$ such that 
$P\in C^\infty (\U_s,\C)$.
\item $P$ has a zero of order at least 3 in $0$.
\item
for each $k\geq 3$ the Taylor's
expansion  of
degree $k$ of $P$ at zero belongs
to $\otimes_{l=3}^k\T_{l}^{\nu,\beta}$.
\end{itemize}
\end{definition}

We now define the class of polynomials in {\bf normal form}:

\begin{definition}\label{normalZ}
Let $k = 2m$ be an even integer. A formal homogeneous polynomial $Z$ of degree $k$ on $\Ps$ is in normal form if it reads
\be Z(z) = \sum_{j \in \bar{\N}^{m}} b_j z_{j_1}z_{-j_1} \dots z_{j_m}z_{-j_m} \ee
i.e. $Z$ depends only on the actions $I_l:=z_lz_{-l}=\xi_l\eta_l$.
\end{definition}
The aim of  the Birkhoff normal form theorem is to reduce a given Hamiltonian of the form $H_0+P$ with $P$ in $\Hs$  to a Hamiltonian of the form $Z+R$ where $Z$ is in normal form and $R$ remains very small, in the sense that it has a zero of high order at the origin.

We now review the properties of polynomials in the class $\T^{\nu,\beta}$.
\begin{proposition}\label{tame}
     Let $k\in \Nb$, $\nu\in[0,+\infty)$, $\beta\in [0,+\infty)$,  
$s\in\R$ with
     $s>\nu+1$,
     and let $P\in \T_{k+1}^{\nu,\beta}$. Then
     \begin{itemize}
     \item[(i)] $P$ extends as a continuous polynomial on $\Ps(\C)$
     and there exists a constant $C>0$ such that for all $z\in \Ps(\C)$
      $$ \va{P(z)}\leq C\norma{z}_{s}^{k+1}$$
     \item[(ii)]  For any  $s'<s+\beta-\frac{1}{2}$,  the
     Hamiltonian vector field $X_{P}$ extends as
     a bounded function from $\Ps(\C)$ to $\mathcal P_{s'}(\C)$.
     Furthermore, for any $s_0\in(\nu+1,s]$, there is $C>0$ such that  
for any
$z\in \Ps(\C)$
\be\label{eqtame-}\norma{X_{P}(z)}_{s'}\leq
C\norma{z}_{s}\norma{z}_{s_{0}}^{(k-1)}.\ee
\item[(iii)] Assume moreover that  $P\in \T_{k+1}^{\nu,\beta,+}$ with  
$\beta>0$, then the
     Hamiltonian vector field $X_{P}$ extends as
     a bounded function from $\Ps(\C)$ to $\Ps(\C)$.
     Furthermore, for any $s_0\in(\nu+1,s]$, there is $C>0$ such that  
for any
$z\in \Ps(\C)$
\be\label{eqtame}\norma{X_{P}(z)}_{s}\leq
C\norma{z}_{s}\norma{z}_{s_{0}}^{(k-1)}.\ee
\item[(iv)] Assume finally that $P\in \T_{k+1}^{\nu,\beta}$ and $P$ is in {\em normal form} in the sense of Definition \ref{normalZ}, then the
     Hamiltonian vector field $X_{P}$ extends as
     a bounded function from $\Ps(\C)$ to $\Ps(\C)$.
     Furthermore, for any $s_0\in(\nu,s]$, there is $C>0$ such that  
for any
$z\in \Ps(\C)$
\be\norma{X_{P}(z)}_{s}\leq
C\norma{z}_{s}\norma{z}_{s_{0}}^{(k-1)}.\ee
     \end{itemize}
\end{proposition}
\begin{remark} 
Since homogeneous polynomials are their own Taylor expansion at $0$, assertions (iii) and (iv) imply that every element of $\T_{k+1}^{\nu,\beta,+}$, and every element of $\T_{k+1}^{\nu,\beta}$ in normal form is in $\Hs$.
\end{remark}

\proof  (i) Let $P$ be an
homogeneous polynomial of degree $k+1$ in $\T_{k+1}^{\nu,\beta}$ and   
 for $z\in \Ps(\C)$ write
\be \label{Pp}
P(z)=\sum_{j\in \Zb^{k+1}}a_{j}\, z_{j_{1}}\ldots z_{j_{k+1}}\ .
\ee
One has, using first \eqref{a} and then that $A(j)\leq 1$, $C(j)\geq 1$,
\begin{eqnarray*}
\va{P(z)}&\leq& C
\sum_{j\in
\Zb^{k+1}}\frac{\mu(j)^{\nu}}{C(j)^\beta}{A(j)^N}\
\prod_{i=1}^{k+1}\va{z_{j_{i}}} \\
&\leq&C
\sum_{j\in
\Zb^{k+1}}\frac{\mu(j)^\nu}{\prod_{i=1}^{k+1}\va{j_{i}}^s}\
\prod_{i=1}^{k+1}\va{j_{i}}^s\va{z_{j_{i}}} \\
&\leq&C
\sum_{j\in
\Zb^{k+1}}\frac{1}{\prod_{i=1}^{k+1}\va{j_{i}}^{s-\nu}}\
\prod_{i=1}^{k+1}\va{j_{i}}^s\va{z_{j_{i}}} \\
&\leq&C
\left( \sum_{l\in
\Zb}\frac{1}{\va{l}^{2s-2\nu}}\right)^{\frac{k+1}{2}}
\norma{z}_{s}^{k+1}
\end{eqnarray*}
where, in the last inequality, we used $k+1$ times the Cauchy-Schwarz
inequality.
Since $s>\nu+1/2$, the last sum converges and the first assertion is  
proved.\\

(ii) We have to estimate the derivative of polynomial $P$ with respect  
to any of its variables. Because of \eqref{a}, given any $N$, we get
\begin{equation*}
\big| \frac{\partial P}{\partial z_l}\big| \leq C_N (k+1) \sum_{j \in  
\Zb^k} \frac{\mu(j,l)^\nu}{C(j,l)^\beta}A(j,l)^N  
|z_{j_1}|\dots|z_{j_k}|\,,
\end{equation*}
where the quantities $\mu(j,l)$, $C(j,l)$ and $A(j,l)$ are computed for  
the $k+1$-tuple made of $j_1, \dots, j_k, l$. Furthermore
\begin{eqnarray}
\nonumber ||X_P(z)||_{s'}^2 &\leq & C \sum_{l \in \Zb} \left( \sum_{j  
\in \Zb^k} \frac{ |l|^{s'} \mu(j,l)^\nu}{C(j,l)^\beta}A(j,l)^N  
|z_{j_1}|\dots|z_{j_k}|\right)^2\\ \nonumber
&\leq& C (k!)^2 \sum_{l \in \Zb} \left( \sum_{j \in \Zb^k_>} \frac{  
|l|^{s'} \mu(j,l)^\nu}{C(j,l)^\beta}A(j,l)^N  
|z_{j_1}|\dots|z_{j_k}|\right)^2\\ \label{XP}
&\leq& C' ||z||_{s_0}^{2(k-3)}\sum_{l \in \Zb} \left( \sum_{|j_1|\geq  
|j_2| \geq |j_3|} \frac{ |l|^{s'} \mu(j,l)^\nu}{C(j,l)^\beta}A(j,l)^N  
|z_{j_1}||z_{j_2}||z_{j_3}|\right)^2\,,
\end{eqnarray}
where $\Zb^k_>$ denotes the set of {\em ordered} $k$-uples $(j_1,  
\dots,j_k)$ such that $|j_1|\geq|j_2|\geq \dots \geq |j_k|$. We used  the following result in the last inequality:
\begin{lemma}\label{lems0}
Given any $s \geq 0$, $s_0 > \frac 12$ and $z \in \ell^2_{s+s_0}$ we  
have
$$ \sum_{j \in \Zb} |j|^s |z_j| \leq C_{s_0} ||z||_{s+s_0}\,.$$
\end{lemma}
\proof
This result is a simple consequence of Cauchy-Schwarz inequality :
$$  \sum_{j \in \Zb} |j|^s |z_j| = \sum_{j \in \Zb}  
\frac{1}{|j|^{s_0}}|j|^{s+s_0} |z_j| \leq \left(\sum_{j \in \Zb}  
\frac{1}{|j|^{2s_0}}\right)^{\frac 12}||z||_{s+s_0}\,.$$ \endproof

Before continuing with the proof  of assertion (ii) of Proposition \ref{tame}, we give  two technical lemmas which give an estimate of $A(j,l)$.

\begin{lemma}\label{lajl}
Given any ordered $k$-tuple $j\in \Zb^k_>$ and $l \in \Zb$, we have $$  
|l| A(j,l) \leq 2 |j_1|\,.$$
\end{lemma}
\proof
It is straightforward if $|l| \leq 2 |j_1|$, since $A(j,l)\leq 1$. If  
not, the order is the following : $|l|> 2|j_1| > |j_1| \geq |j_2|$ and
$$ |l| A(j,l) = \frac{|l| \sqrt{|j_1j_2|}}{\sqrt{|j_1j_2|} + |l| - |j_1|}  
\leq \frac{|l| \sqrt{|j_1j_2|}}{ |l|/2} \leq 2 |j_1|\,,$$
and the lemma is proved.
\endproof
\begin{lemma}\label{ajlmax}
Given any ordered $k$-uple $j\in \Zb^k_>$ and $l \in \Zb$ we have
$$ A(j,l) \leq \tilde{A}(j_1,j_2,l) := \left\{ \begin{array}{ccc} 2  
\frac{|j_2|}{|l|+|j_1|-|j_2|} & \mbox{if} & |l|\leq |j_2| \,,\\
2 \frac{\sqrt{|lj_2|}}{\sqrt{|lj_2|} + ||j_1|-|l||} & \mbox{if} & |l|  
\geq |j_2|\,.\end{array}\right.$$
\end{lemma}
\proof If $|l| > 2 |j_1|$, $A(j,l)$ reads :
$$ A(j,l) = \frac{\sqrt{|j_1j_2|}}{\sqrt{|j_1j_2|}+|l|-|j_1|}\,.$$
We can write :
\begin{eqnarray*}
\sqrt{|j_1j_2|}+|l|-|j_1| & = & \sqrt{|lj_2|}+|l|-|j_1| -  
\sqrt{|j_2|}(\sqrt{|l|}-\sqrt{|j_1|})\\
& =  & \sqrt{|lj_2|}+|l|-|j_1| - \sqrt{|j_2|}\frac{|l| -  
|j_1|}{\sqrt{|l|}+\sqrt{|j_1|}}\\
& \geq & \sqrt{|lj_2|}+|l|-|j_1| -\sqrt{|j_1|}\frac{|l| -  
|j_1|}{\sqrt{|l|}+\sqrt{|j_1|}}\\
& \geq &\sqrt{|lj_2|} + \frac{\sqrt{2}}{\sqrt{2}+1}(|l|-|j_1|)
\end{eqnarray*}
Hence,
$$ A(j,l) \leq  
\frac{1+\sqrt{2}}{\sqrt{2}}\frac{\sqrt{|j_1j_2|}}{\sqrt{|j_1j_2|}+|l|- 
|j_1|} \leq 2 \frac{\sqrt{|lj_2|}}{\sqrt{|lj_2|}+|l|-|j_1|}\,.$$

If $|j_2| \leq |l| \leq 2 |j_1|$, then $B(j,l)^2 = |j_2|  
\min(|l|,|j_1|) \in \left[ \frac{|l j_2|}{2}, |lj_2| \right]$,  
therefore
$$ A(j,l) \leq \frac{\sqrt{|lj_2|}}{1/\sqrt{2} \sqrt{|lj_2|} +  
||l|-|j_1||} \leq 2 \frac{\sqrt{|lj_2|}}{\sqrt{|lj_2|} +  
||l|-|j_1||}\,.$$
Finally, if $|l| \leq |j_2|$ we get
$$ A(j,l) = \frac{\sqrt{|lj_2|}}{\sqrt{|lj_2|} + |j_1|-|j_2|} \leq 2  
\frac{|j_2|}{|l|+|j_1|-|j_2|}\,,$$
and this ends the proof of Lemma \ref{ajlmax}.
\endproof

To continue with the proof of assertion (ii) of Proposition \ref{tame},  
we define  $0<\varepsilon < s -s'-\frac 12$, and $N = s + 1  
+\varepsilon$. In view of \eqref{XP}, we may decompose :
\begin{equation}
  || X_P(z)||_{s'}^2  \leq C \sum_{l \in \Zb} \left( T_1(l) +T_2(l)  
\right)^2\,,
\end{equation}
with
\begin{eqnarray*}
T_1(l) & = & \sum_{|j_1|\geq |j_2|\geq|j_3|, |j_2| > |l|}  
\frac{|l|^{s'}|\mu(j,l)|^\nu}{\max(|j_1|,|l|)^{\beta}}A(j,l)^{N}|z_{j_1 
}| |z_{j_2}||z_{j_3}|\\
T_2(l) & = & \sum_{|j_1|\geq |j_2|\geq|j_3|, |j_2| \leq |l|}  
\frac{|l|^{s'}|\mu(j,l)|^\nu}{\max(|j_1|,|l|)^{\beta}}A(j,l)^{N}|z_{j_1 
}| |z_{j_2}||z_{j_3}|\,.
\end{eqnarray*}
Since $A(j,l) \leq 1$ and $N > \frac12 + s' + \varepsilon$, we may  
estimate $T_1(l)$ using Lemmas \ref{lajl} and \ref{ajlmax} :
\begin{eqnarray*}
T_1(l) & \leq & C \sum_{|j_1|\geq |j_2|\geq|j_3|, |j_2| > |l|}  
|j_1|^{s'}|j_2|^\nu \tilde{A}(j_1,j_2,l)^{\frac 12 + \varepsilon  
}|z_{j_1}| |z_{j_2}||z_{j_3}|\\
& \leq & C ||z||_{s_0} \sum_{|j_1|\geq |j_2|, |j_2| > |l|}  
\frac{1}{|l|^{\frac 12 + \varepsilon}}|j_1|^{s-\frac 12 -  
\varepsilon}|z_{j_1}| |j_2|^{\nu + \frac 12 + \varepsilon}|z_{j_2}|\\
& \leq & C ||z||_{s_0}\frac{1}{|l|^{\frac 12 + \varepsilon}}||z||_s  
||z||_{\nu+ 1 + 2 \varepsilon}\,,
\end{eqnarray*}
hence $T_1(l)$ is an $\ell^2$-sequence, whose $\ell^2$-norm is bounded  
above by $C||z||_{s_0}^2 ||z||_s$ if we assume that $s_0 > \nu + 1 + 2\varepsilon$.  
Concerning $T_2(l)$, using Lemmas \ref{lajl} and \ref{ajlmax}, we  
obtain
\begin{eqnarray*}
T_2(l) & \leq & C \sum_{|j_1|\geq |j_2|\geq|j_3|, |j_2| \leq |l|}  
\frac{1}{|l|^{s-s'+\beta}}|j_1|^s |j_2|^\nu  
\tilde{A}(j_1,j_2,l)^{N-s}|z_{j_1}| |z_{j_2}||z_{j_3}|\\
& \leq & C \frac{||z||_{s_0}}{|l|^{s-s'+\beta}} \sum_{|j_1|\geq |j_2|,  
|j_2| \leq |l|} \left( \frac{\sqrt{|lj_2|}}{1 + ||j_1|-|l||}\right)^{1  
+ \varepsilon}|j_1|^s |z_{j_1}| |j_2|^\nu |z_{j_2}|\\
& \leq & C  \frac{||z||_{s_0}}{|l|^{s-s'+\beta- (1 +  
\varepsilon)/2}}\left(\sum_{j_2 \in \Zb} |j_2|^{\nu + (1 + \varepsilon)/2}  
|z_{j_2}| \right)\sum_{j_1 \in \Zb} \frac{|j_1|^s |z_{j_1}|}{(1 + | |j_1| -  
|l| | )^{1+\varepsilon}}.
\end{eqnarray*}
The last sum in $j_1$ is a convolution product of the $\ell^2$-sequence  
$|j_1|^s |z_{j_1}|$ and the $\ell^1$-sequence $\frac{1}{(1 + |j_1|  
)^{1+\varepsilon}}$ and thus a $\ell^2$-sequence with respect to the index $l$, whose $\ell^2$-norm is bounded by $||z||_s$ . Choosing $\varepsilon>0$ in such a way that  
$s-s'+\beta- (1 + \varepsilon)/2>0$, the sequence  $T_2(l)$ is in  
$\ell^2$, with a norm bounded by
\begin{equation*}
||T_2|| \leq C ||z||_{s_0}||z||_{\nu+(1+\varepsilon)/2} ||z||_s \leq C  
||z||_{s_0}^2||z||_s\,,
\end{equation*}
with $s_0 > \nu + (1+\varepsilon)/2$. Collecting the estimates for  
$T_1$ and $T_2$, we obtain the desired inequality.\\

(iii) We define $0<\varepsilon <  1/12$ and $N = s +\frac 12  
+\varepsilon$. We have, as in (ii), this first estimate
\begin{equation}
||X_P(z)||_{s}^2 \leq  C ||z||_{s_0}^{2(k-3)}\sum_{l \in \Zb} \left(  
\sum_{|j_1|\geq |j_2| \geq |j_3|} \frac{  
|l|^{s}\mu(j,l)^\nu}{C(j,l)^\beta(1+S(j,l))}A(j,l)^{N}  
|z_{j_1}||z_{j_2}||z_{j_3}|\right)^2\,.
\end{equation}
As in (ii), we may also decompose the sum on $j_1$, $j_2$ and $j_3$ into two  
pieces, $T_1^+(l)$ collecting all the terms with $|j_2|>|l|$ and  
$T_2^+(l)$ collecting those with $|j_2| \leq |l|$. Following (ii),  
since $C(j,l)\geq 1$ and $1 + S(j,l)\geq 1$, we obtain for $T_1^+$ :
\begin{eqnarray*}
T_1^+(l) & \leq & C \sum_{|j_1|\geq |j_2|\geq|j_3|, |j_2| > |l|}  
|l|^{s}|j_2|^\nu A(j,l)^{N}|z_{j_1}| |z_{j_2}||z_{j_3}|\\
& \leq & C \sum_{|j_1|\geq |j_2|\geq|j_3|, |j_2| > |l|}  
|l|^{1/2+\varepsilon}|j_1|^{s-1/2-\varepsilon}|j_2|^\nu  
A(j,l)^{N-(s-1/2-\varepsilon)}|z_{j_1}| |z_{j_2}||z_{j_3}|\\
& \leq & C ||z||_{s_0} \sum_{|j_1|\geq |j_2|, |j_2| >  
|l|}|l|^{1/2+\varepsilon}|j_1|^{s-1/2-\varepsilon}|j_2|^\nu  
\tilde{A}(j_1,j_2,l)^{N-(s-1/2-\varepsilon)}|z_{j_1}| |z_{j_2}|\\
& \leq & C ||z||_{s_0} \sum_{|j_1|\geq |j_2|, |j_2| >  
|l|}|l|^{s-N}|j_1|^{s-1/2-\varepsilon}|j_2|^{\nu+N-(s-1/2- 
\varepsilon)}|z_{j_1}| |z_{j_2}|\\
& \leq & C ||z||_{s_0} \frac{1}{|l|^{\frac 12 +  
\varepsilon}}\sum_{|j_1|\geq |j_2|, |j_2| >  
|l|}|j_1|^{s-1/2-\varepsilon}|j_2|^{\nu+N-(s-1/2-\varepsilon)}|z_{j_1}|  
|z_{j_2}|\\
& \leq & C ||z||_{s_0}\frac{1}{|l|^{\frac 12 + \varepsilon}}||z||_s  
||z||_{\nu+ 1 + 2 \varepsilon}\,,
\end{eqnarray*}
hence $T_1^+(l)$ is a $\ell^2$-sequence, whose $\ell^2$-norm is bounded  
above by $C||z||_{s_0}^2 ||z||_s$ if $s_0 > \nu + 1 + 2\varepsilon$.

The estimate on $T_2^+$ will need all factors assigned in the  
definition of $\T^{\nu,\beta,+}$:
\begin{eqnarray*}
T_2^+(l) &\leq & C \sum_{|j_1|\geq |j_2|\geq|j_3|, |j_2| \leq |l|}  
\frac{|j_1|^s |j_2|^\nu}{\max(j_1,l)^{\beta}(1 + ||j_1|-|l||)}  
\tilde{A}(j_1,j_2,l)^{N-s}|z_{j_1}| |z_{j_2}||z_{j_3}|\\
& \leq & C ||z||_{s_0} \sum_{|j_1|\geq |j_2|, |j_2| \leq |l|} \left(  
\frac{\sqrt{|lj_2|}}{1 +  
||j_1|-|l||}\right)^{\varepsilon}\frac{1}{|l|^{\beta}(1 + ||j_1|-|l||)}  
|j_1|^s |z_{j_1}| |j_2|^\nu |z_{j_2}|\\
& \leq & C ||z||_{s_0} \frac{1}{|l|^{\beta- \varepsilon/2}}\sum_{j_2  
\in \Zb} |j_2|^{\nu + \varepsilon/2} |z_{j_2}| \sum_{j_1 \in \Zb}  
\frac{|j_1|^s |z_{j_1}|}{(1 + | |j_1| - |l| | )^{1+\varepsilon}}\,.
\end{eqnarray*}
Once again, the last sum in $j_1$ is a convolution product of the  
$\ell^2$ sequence $|j_1|^s |z_{j_1}|$ and the $\ell^1$ sequence  
$\frac{1}{(1 + |j_1| )^{1+\varepsilon}}$. Choosing $\varepsilon>0$ in  
such a way that $\beta- \varepsilon/2>0$, the sequence  $T_2^+(l)$ is  
in $\ell^2$, with a norm bounded by
\begin{equation*}
||T_2|| \leq C ||z||_{s_0}||z||_{\nu+(1+\varepsilon)/2} ||z||_s \leq C  
||z||_{s_0}^2||z||_s\,,
\end{equation*}
with $s_0 > \nu + (1+\varepsilon)/2$. Collecting the estimates for  
$T_1^+$ and $T_2^+$, we obtain the stated inequality.\\

(iv) Let $k+1=2m$. As in (ii), we obtain 
\be\label{aqw1} || X_P||_s^2 \leq C \sum_{l \in \Zb} \left( \sum_{j \in \Nb_>^{m-1}} |l|^s |z_l| \frac{ \mu(j,j,l,l)^\nu}{C(j,j,l,l)^\beta}|z_{j_1}||z_{-j_1}|\dots |z_{j_{m - 1}}| |z_{-j_{m - 1}}| \right)^2,\ee
using the same convention for $\mu(j,j,l,l)$ and $C(j,j,l,l)$ as for $\mu(j,l)$ and $C(j,l)$: as an example,  $\mu(j,j,l,l)$ is the third biggest integer among \\ $|j_1|, |j_1|, \dots |j_{m-1}|, |j_{m-1}|, |l|$ and $|l|$, that is, if $j$ is ordered, either $\mu(j,j,l,l)=|j_1|$, and in this case $C(j,j,l,l)=\val$,  or $\mu(j,j,l,l)=\val$ and in this case $C(j,j,l,l)=|j_1|$. Notice that $A(j,j,l,l)=1$ does not help for this computation. The sum over $j$ can be decomposed into two parts :
\begin{eqnarray*} \lefteqn{\sum_{j \in \Nb_>^{m-1}, j_1 \leq l} |l|^s |z_l| \frac{ \mu(j,j,l,l)^\nu}{C(j,j,l,l)^\beta}|z_{j_1}||z_{-j_1}|\dots |z_{j_{m - 1}}| |z_{-j_{m - 1}}|} \\& \leq & \sum_{j \in \Nb_>^{m-1}, j_1 \leq l} |l|^s |z_l| \frac{|j_1|^\nu}{|l|^\beta} |z_{j_1}||z_{-j_1}|\dots |z_{j_{m - 1}}| |z_{-j_{m - 1}}|\\
& \leq  & |l|^{s-\beta} |z_l| \sum_{j_1} j_1^\nu |z_{j_1}||z_{-j_1}|  ||z||_0^{2(m-2)}\\
& \leq & |l|^{s-\beta} |z_l| ||z||_{\nu/2}^2 ||z||_0^{2(m-2)},
\end{eqnarray*}
and
\begin{eqnarray*} \lefteqn{\sum_{j \in \Nb_>^{m-1}, j_1 > l} |l|^s |z_l| \frac{ \mu(j,j,l,l)^\nu}{C(j,j,l,l)^\beta}|z_{j_1}||z_{-j_1}|\dots |z_{j_{m - 1}}| |z_{-j_{m - 1}}|}\\ & \leq & |l|^s |z_l| \sum_{j \in \Nb_>^{m-1}} |j_1|^{\nu - \beta} |z_{j_1}||z_{-j_1}|\dots |z_{j_{m - 1}}| |z_{-j_{m - 1}}|\\ & \leq & |l|^s |z_l| ||z||_{(\nu-\beta)/2}^2 ||z||_0^{2(m-2)}
\end{eqnarray*}
Inserting these two estimates  in (\ref{aqw1}) we get \eqref{eqtame}.
\endproof

The second essential property satisfied by polynomials in $\T_{k}^{\nu,\beta}$
is captured in the following
\begin{proposition}\label{bra} Let $k_1,k_2\geq2$, $\nu_1,\nu_2\geq0$  
and $\beta>0$
    The map $
  (P,Q)\mapsto \{P,Q\}$ defines a continuous map from
  $\T_{k_{1}+1}^{\nu_{1},\beta,+}\times\T_{k_{2}+1}^{\nu_{2},\beta}$ to $
  \T_{k_{1}+k_{2}}^{\nu',\beta}$ for
   $\nu'=2(\nu_{1}+\nu_{2})+1$.
\end{proposition}
\proof
We assume that
$P\in \T_{k_{1}+1}^{\nu_{1},\beta,+}$ and $Q\in  
\T_{k_{2}+1}^{\nu_{2},\beta}$
are homogeneous polynomials and we write
$$
P(z)=\sum_{j\in \Zb^{k_{1}+1}}a_{j}\, z_{j_{1}}\ldots z_{j_{k_{1}+1}}
$$
  and
$$
Q(z)=\sum_{i\in \Zb^{k_{2}+1}}b_{i}\, z_{i_{1}}\ldots z_{i_{k_{2}+1}}\ .
$$
In view of the symmetry of the   estimate \eqref{a} with respect to
the involved indices, one easily obtains
$$\{P,Q\}(z)=\sum_{(j,i)\in \Zb^{k_{1}+k_{2}}}c_{j,i}\, z_{j_{1}}\ldots
z_{j_{k_{1}}}z_{i_{1}}\ldots z_{i_{k_{2}}}$$
with
$$
\va{c_{j,i}}\leq c_{N,N'}\sum_{l\in \Zb}
\frac{\mu(j,l)^{\nu_{1}}}{C(j,l)^\beta(1+S(j,l))}A(j,l)^{N}
\frac{\mu(i,l)^{\nu_{2}}}{C(i,l)^\beta}A(i,l)^{N'}.
$$
Therefore it remains to prove that, for each $M\geq 1$, there exist  
$N,N'\geq1$,  $C>0$ such that for
all  $j\in \Zb^{k_{1}}$ and all $i\in \Zb^{k_{2}}$,
\be \label{+}
\sum_{l\in \Zb}
\frac{\mu(j,l)^{\nu_{1}}}{C(j,l)^\beta(1+S(j,l))}A(j,l)^{N}
\frac{\mu(i,l)^{\nu_{2}}}{C(i,l)^\beta}A(i,l)^{N'}
\leq C \frac{\mu(j,i)^{\nu'}}{C(i,j)^\beta}{A(j,i)^{M}}
\ee
with $\nu'=2(\nu_{1}+\nu_{2})+1$.\\
In order to simplify the notations, and because it does not change the  
estimates
of \eqref{+},
we will assume that $k_{1}=k_{2}=k$. We can also assume by symmetry that
\begin{itemize}
\item all the indices are positive:
$j_{1},\ldots,j_{k},i_{1},\ldots,i_{k}\geq 1$.
\item $j$ and $i$ are ordered:
$j_{1}\geq \ldots \geq j_{k}$ and $i_{1}\geq \ldots \geq i_{k}$.
\end{itemize}
We begin with two technical lemmas whose proofs are postponed at the end  
of this proof.
\begin{lemma}\label{11}
There is a constant $C>0$ such that for any  $j\in \Zb^{k_{1}}$, $i\in  
\Zb^{k_{2}}$ and $l\in \Zb$ we have
\be\label{estim11}
A(j,l)^2A(i,l)^2\leq C A(i,j).
\ee
\end{lemma}
\begin{lemma}\label{12}
There is a constant $C>0$ such that for any  $j\in \Zb^{k_{1}}$, $i\in  
\Zb^{k_{2}}$ and $l\in \Zb$ we have
\be\label{estim12}
\max(\mu(j,l)A(i,l)^2,\mu(i,l)A(j,l)^2)\leq C \mu(i,j)^2.
\ee
\end{lemma}

Using these lemmas, in order to prove \eqref{+}, it suffices to prove
$$
\sum_{l\in \Zb}
\frac{1}{C(j,l)^\beta(1+S(j,l))}
\frac{A(i,l)^2}{C(i,l)^\beta}
\leq C \frac{\mu(j,i)}{C(i,j)^\beta}.
$$
Noticing that $C(i,l)C(j,l)\geq C(i,j)l$, it suffices to verify that
$$
\sum_{l\in \Zb}
\frac{A(i,l)^2}{(1+S(j,l))l^{\beta}}
\leq C {\mu(j,i)}.
$$
Decompose the sum in two parts, $I_1=\sum_{l>j_2}$ and  
$I_2=\sum_{l\leq j_2}$. For the first sum we have
$$
I_1=\sum_{l>j_2}
\frac{A(i,l)^2}{(1+S(j,l))l^{\beta}}
\leq \sum_{l\in\Zb}
\frac{1}{(1+\va{l-j_1})l^{\beta}}\leq C,
$$
while for the second one
$$
I_2=\sum_{l\leq j_2}\frac{A(i,l)^2}{(1+S(j,l))l^{\beta}}\leq
\sum_{l\leq j_2}\frac{A(i,l)^2}{l^{\beta}}.
$$
In this last sum, if $j_2<\mu(i,j)$, then
$$
I_{2}\leq
j_2\leq \mu(i,j).
$$
On the other hand, if $\mu(i,j)\leq j_2$, then we decompose the $I_2$ sum  
in two parts,
$I_{2,1}=\sum_{l<2i_1}$ and $I_{2,2}=\sum_{l\geq 2i_1}$. Since $i_1\leq  
\mu(i,j)=\max(i_1,j_3)$ we have
$$I_{2,1}=\sum_{l\leq 2i_1}\frac{A(i,l)^2}{l^{\beta}}\leq 2i_1\leq  
2\mu(i,j).$$
Finally, when $l\geq 2i_1$ we have $S(i,l)\geq l/2$ and  
$B(i,l)^2=i_1i_2\leq i_2l/2\leq \mu(i,j)l/2$ and thus $A(i,l)\leq  
\sqrt{2\mu(i,j)}l^{-1/2}$ which leads to
$$
I_{2,2}=\sum_{2i_1\leq l\leq j_2}\frac{A(i,l)^2}{l^{\beta}}\leq  
C\mu(i,j)\sum_{l\in \Nb}
\frac 1 {l^{1+\beta}}
\leq C\mu(i,j).
$$
\endproof

{\it Proof of lemma \ref{11} ---}
The estimate \eqref{estim11}, being symmetric with respect to $i$ and  
$j$, we can assume that $j_1\geq i_1$.
  We consider three cases, depending of the position of $l$ with respect  
to $i_1$ and $j_1$.\\
  {\bf First case $l\geq j_1$:}\\
  We have $S(i,l)=\va{i_1-l}\geq \va{i_1-j_1}\geq S(i,j)$ and  
$B(i,l)=(i_1i_2)^{1/2}\leq B(i,j)$. Therefore
  $$A(i,l)=\frac{B(i,l)}{B(i,l)+S(i,l)}\leq  
\frac{B(i,l)}{B(i,l)+S(i,j)}\leq \frac{B(i,j)}{B(i,j)+S(i,j)}=A(i,j),$$
and using $A(j,l)\leq 1$, \eqref{estim11} is proved.

\medskip
\noindent {\bf Second case} $l\leq i_1$:\\
  Similarly as in the first case, we have $S(j,l)\geq S(i,j)$ and  
$B(j,l)=(j_2\max(j_3,l))^{1/2}\leq (j_2\max(j_3,i_1))^{1/2}\leq B(i,j)$  
and thus
  $$A(j,l)\leq A(i,j).$$
  {\bf{Third case}} $i_1<l<j_1$:\\
  That is the most complicated case and we have to distinguish whether
  $i_1\geq j_2$ or not.\\
{\it Subcase 1. } $i_1\geq j_2$:\\
We have $B(i,l)\leq B(i,j)$ thus if $S(i,l)=\va{i_1-l}\geq\frac 1 2  
\va{i_1-j_1}=\frac 1 2 S(i,j)$ we obtain $A(j,l)\leq2A(i,j)$ and  
\eqref{estim11} holds true. Now if $S(i,l)<\frac 1 2 S(i,j)$ then  
$S(j,l)\geq \frac 1 2 S(i,j)$
since $S(i,l)+S(j,l)\geq S(i,j)$. Furthermore, if $B(j,l)\leq B(i,j)$  
then
$$A(j,l)=\frac{B(j,l)}{B(j,l)+S(j,l)}\leq 2  
\frac{B(j,l)}{B(j,l)+S(i,j)}\leq 2 \frac{B(i,j)}{B(i,j)+S(i,j)}=2  
A(i,j)\,,
$$
and \eqref{estim11} holds. If $B(j,l)> B(i,j)$, then using
$$B(j,l)^2=j_2l =j_2i_1 +j_2(l-i_1)\leq B(i,j)^2+j_2S(i,l)\leq  
B(i,j)^2+\frac 1 2 B(i,j)S(i,j)\,,$$
we deduce
$$A(j,l)^2 \leq \frac{B(j,l)^2}{(B(i,j) + \frac 12 S(i,j))^2} \leq  
2\frac{B(i,j)^2+ B(i,j)S(i,j)}{(B(i,j)+S(i,j))^2}\leq  
2(A(i,j)^2+A(i,j))\,,$$
thus \eqref{estim11} is also satisfied in this case, since $A(i,j) \leq  
1$.

\medskip

\noindent {\it Subcase 2. } $i_1\leq j_2$:\\
We still have $B(i,l)\leq B(i,j)$ thus if furthermore $S(i,j)\leq  
2S(i,l)$ then $A(i,l)\leq 2A(i,j)$ and \eqref{estim11} is true. So we  
assume $2S(i,l)<S(i,j)$ which implies $S(i,j)\leq 2S(j,l)$ since  
$S(i,l)+S(j,l)\geq S(i,j)$. If furthermore $l\leq j_3$,  
$B(j,l)=B(j)\leq B(i,j)$ and thus $A(j,l)\leq 2A(i,j)$ and  
\eqref{estim11} is again true. So we assume $j_3\leq l$ and we have
$$B(j,l)^2=lj_2=i_1j_2+j_2(l-i_1)\leq B(i,j)^2 +j_2S(i,l).$$
If $S(i,l)\leq l/2$ then we deduce $B(j,l)^2\leq 2B(j,l)^2$ and  
\eqref{estim11} is satisfied. It remains to consider the case  
$S(i,l)>l/2$ which implies $i_1<l/2$ and thus
\be\label{A1}
A(i,l)\leq \frac{i_1}{i_1+l/2}\leq 2\frac{i_1}l .
\ee
Let $n\geq 1$ such that $\frac{l}{2^{n+1}}\leq i_1\leq \frac{l}{2^{n}}$  
we get from \eqref{A1}
\be\label{A2}
A(i,l)\leq \frac{1}{2^{n-1}} .
\ee
On the other hand
\be \label{A3}
A(j,l)\leq 2 \frac{(lj_2)^{1/2}}{(lj_2)^{1/2}+S(i,j)}\leq 2  
\frac{(lj_2)^{1/2}}{(i_1j_2)^{1/2}+S(i,j)}
\ee
and
\be \label{A4}
A(i,j)\geq  \frac{(i_1j_2)^{1/2}}{(i_1j_2)^{1/2}+S(i,j)}\geq \frac  
1{2^{n+1}}\frac{(lj_2)^{1/2}}{(i_1j_2)^{1/2}+S(i,j)}.
\ee
Combining \eqref{A2}, \eqref{A3} and \eqref{A4} we conclude
$$A(i,l)A(j,l)\leq 8 A(i,j).$$
\endproof
\noindent {\it Proof of lemma \ref{12} ---}
The estimate \eqref{estim12} being symmetric with respect to $i$ and  
$j$, we can assume $j_1\geq i_1$. If furthermore $i_1\geq j_2$ then one  
easily verifies that
$$\mu(i,l)\leq \mu(i,j) \mbox{ and } \mu(j,l)\leq \mu(i,j)$$ and  
estimates \eqref{estim12} is satisfied.\\
In the case $j_1\geq j_2 \geq i_1$ we still have $\mu(i,l)\leq \mu(i,j)  
$ but $\mu(j,l)$ could be larger than $\mu(i,j)$. Actually if  $  
\mu(j,l)\leq 2\mu(i,j)$, estimates \eqref{estim12} is still trivially  
satisfied. Therefore it remains to consider the case where $ \mu(j,l)>  
2\mu(i,j)$. Remark that in this case $i_1\leq\mu(i,j)\leq  
\frac{\mu(j,l)} 2 \leq l/2$ and thus $S(i,l)=\va{i_1 -l}\geq l/2$   
which  leads to
$$A(i,l)\leq \frac{(i_1i_2)^{1/2}}{S(i,l)}\leq\frac{  
(2i_2)^{1/2}}{l^{1/2}}\leq \frac{ (2\mu(i,j))^{1/2}}{l^{1/2}}.$$
Using this last estimates one gets
$$\mu(j,l)A(i,l)^2\leq l A(i,l)^2\leq 2\mu(i,j)^2.$$
\endproof

We end this section with a proposition concerning Lie transforms of homogeneous polynomials $\chi \in \T^{\delta,\beta,+}_l$, i.e. time 1 flow of the Hamiltonian vector field  $X_\chi$.
\begin{proposition}\label{compo}
Let $\chi$ be a real homogeneous polynomial in $\T_{l}^{\delta,\beta,+}$
with $\delta \geq 0$, $\beta>0$, $l\geq 3$   take $s> s_{1}:=\delta +3/2$ and
denote by $\phi$ the  Lie transform associated with $\chi$. We have
\begin{itemize}
\item[(i)]  $\phi$ is an analytic canonical  
transformation from an open ball $B_\epsilon$ of center $0$ and radius  
$\epsilon$   in  $\Ps$ into  the open ball $B_{2\epsilon}$ in $\Ps$  
satisfying
\begin{equation}\label{estim:phi}\norma{\phi(z)-z}_s\leq C_s \norma{z}^2_s \mbox{  
for any } z\in B_\e.\end{equation}
In particular if $F\in \Hs$ with $s>s_1$
  then
$F\circ \phi\in\Hs$. Furthermore, if $F$ is real then $F\circ \phi$ is real too.
\item[(ii)] Let $P\in \T^{\nu,\beta}_{n} \cap \Hs$, $\nu \geq 0$,
$n\geq 3$ and fix $r\geq n$ an integer. Then
$$P\circ \phi = Q_{r}+ R_{r}$$
where: \\
- $Q_{r}$ is a polynomial of degree at most $r$, belonging to  
$\T^{\nu',\beta}\cap \Hs$
with\\ $\nu'= 2^{r-n} \nu +(2^{r-n}-1)(2\beta +1)$,\\
- $R_{r}$ is a  Hamiltonian in $\T^{\nu'',\beta}\cap \Hs$ with $\nu''= 2^{r-n+1}\nu +(2^{r-n+1}-1)(2\beta +1)$, having a zero of  
order $r+1$
at the origin.
\end{itemize}
\end{proposition}
\proof
(i)  Since $\chi\in \T_{l}^{\delta,\beta,+}$, by Proposition
\ref{tame}(iii), $X_\chi\in C^\infty (\Ps,  
\Ps)$ for $s> s_{1}=\delta +3/2$.  
In
particular, for $s>s_{1}$, the flow $\Phi^t$ generated by the  
vector field $X_{\chi}$
transports an open neighborhood  of the origin in $\Ps$ into an open  
neighborhood of the origin in $\Ps$. Notice that since $\chi$ is real, $\Phi^t$ transports the "real part" of $\Ps$, $\{(\xi,\bar\xi)\in \Ps\}$, into itself.
Furthermore one has for $z\in \Ps$ small enough
$$
\Phi^t(z)-z= \int_0^t X_{\chi}(\Phi^{t'}(z))dt'
$$
and since $\chi$ has a zero of order $3$ at least, one gets by Proposition
\ref{tame}(iii),
$$
\norma{\Phi^t(z)-z}_s\leq C_s  
\int_0^t\norma{\Phi^{t'}_\chi(z)}_s^2dt'.
$$
Then, by a classical continuity argument, there exists $\epsilon>0$  
such that
the flow $B_\e \ni z \mapsto
\Phi_\chi^t(z)\in B_{2\e}$ is well defined and smooth for $0\leq t\leq  
1$. Furthermore, the Lie transform  $\phi=\Phi^1$ satisfies  
\eqref{estim:phi}.\\
On the other hand, by simple composition we get that if $F\in \Hs$ with $s>s_{1}$, then $F \circ \phi \in C^\infty(B_\epsilon, \C)$. In view of the formula
$$X_{F\circ \phi}(z)=(D\phi(z))^{-1} X_{F}(\phi(z)),$$ we deduce
that $X_{F\circ \phi} \in C^\infty (B_\epsilon,  
\Ps)$. We now have to check the properties concerning the Taylor polynomials of $F \circ \phi$. Denoting by $F_k$ (resp. $(F\circ\phi)_k$) the homogeneous polynomial of degree $k$ appearing in the Taylor expansion of $F$ (resp $F\circ\phi$), and putting $F_k^{[0]} = F_k$, $F_k^{[j+1]}= \{F_k, \chi\}$, we have 
$$ (F\circ\phi)_k(z) = \sum_{j\geq 0, k'\geq 0, k'+ j(l-2)= k} F_{k'}^{[j]}(z), $$
since $\chi$ is itself a homogeneous polynomial of degree $l$. It is then sufficient to prove that the Poisson bracket of a homogeneous polynomial $F_k$ in $\Hs$ with $\chi$ stays in $\Hs$. \\
Using the (constant) symplectic form $\omega$ on $\Ps$, we get $\{F_k, \chi\} (z) = \omega(X_{F_k},X_\chi)$, and so $\{F_k, \chi\} \in C^\infty(B_\epsilon, \C)$. Moreover 
$$X_{\{F_k, \chi\}} (z)= [ X_{F_k}, X_\chi ] = \lim_{t \rightarrow 0}\frac 1t ( X_{F_k} - \Phi_*^t(X_{F_k}))(z)\,.$$
Since $\Phi^t$ is the flow of the regular Hamiltonian $\chi \in C^\infty(B_\epsilon, \R)$, the Cauchy Lipschitz theorem implies that the mapping $(t,z) \mapsto \Phi_*^t(X_{F_k})(z)$ is in $C^\infty( [-1,1]\times B_\epsilon, \Ps)$. Now, $X_{\{F_k, \chi\}}$ is nothing else but the time derivative of this mapping at time $0$, hence $X_{\{F_k, \chi\}} \in C^\infty(B_\epsilon, \Ps)$ and the claim is proved.

(ii)
By a direct calculus one has
$$
\frac{d^k}{dt^k}P\circ \Phi^t (z)\big|_{t=0}=P^{[k]}(z)
$$
with the same notation $P^{[k+1]}=\{P^{[k]},\chi \}$ and $P^{[0]}=P$. Therefore
applying the Taylor's formula to $P\circ \Phi^t (z)$ between $t=0$
and $t=1$ we deduce
\be\label{Pphi}
P\circ \phi(z)=
\sum_{k=0}^{r-n}\frac{1}{n!}P^{[k]}(z)+\frac{1}{(r-n)!}\int_{0}^1
(1-t)^r P^{[r-n+1]}(\Phi^t(z))dt.
\ee
Notice that $P^{[k]}(z)$ is a homogeneous
polynomial of degree $n+k(l-2)$ and, by
Proposition \ref{bra}, $P^{[k]}(z)\in \T^{2^k \nu
+(2^k -1)(2\delta+1),\beta}$. Moreover $P^{[k]}(z)$ is a homogeneous polynomial in the Taylor expansion of $P \circ \phi \in \Hs$, hence it is in $\Hs$.
Therefore \eqref{Pphi} decomposes in the sum of a polynomial
of degree $r$ in $\T^{\nu',\beta}_{r}$, and a function in $\Hs$ having a  
zero
of degree $r+1$ at the origin. 
\endproof

\subsection{The Birkhoff normal form theorem}
We start with the resolution of the homological equation and then state the normal form theorem.
\begin{lemma}\label{homo}
Let $\nu \in [0, +\infty)$ and assume that the frequency vector of
$H_{0}$ is strongly non resonant (see Definition \ref{SNR}).
Let $Q$ be a homogeneous real  polynomial of degree $k$ in
$\T_{k}^{ \nu,\beta}$, there exist $\nu'>\nu$, and $Z$ and $\chi$ two   
homogeneous real polynomials of degree $k$, respectively  
in $\T_{k}^{
\nu',\beta}$ and  $\T_{k}^{
\nu',\beta,+}$, which  satisfy
\be \label{homolo}\{H_{0},\chi\}+Q=Z\ee
and 
\be \label{normal}\{Z,I_{j}\}=0\quad \forall j\geq 1\ee
and thus $Z$ is in normal form. 
Furthermore, $Z$ and $\chi$ both belong to $\Hs$ for $s> \nu' +1$.
\end{lemma}
\proof
For $j\in \Nb^{k_{1}}$ and $l\in \Nb^{k_{2}}$ with $k_{1}+k_{2}=k$ we
denote
$$\xi^{(j)}\eta^{(l)}=\xi_{j_{1}}\ldots
\xi_{j_{k_{1}}}\eta_{l_{1}}\ldots\eta_{l_{k_{2}}}.$$
One has
$$
\{H_{0}, \xi^{(j)}\eta^{(l)}
\}
=-i\Omega(j,l) \xi^{(j)}\eta^{(l)}$$
with
$$
\Omega(j,l):=
\omega_{j_{1}}+\ldots+\omega_{j_{k_{1}}}-\omega_{l_{1}}-\ldots
-\omega_{l_{k_{2}}}.$$
Let $Q\in \T^{ \nu,\beta}_{k}$
$$
Q=
\sum_{(j,l)\in \Nb^{k}}a_{jl}\xi^{(j)}\eta^{(l)}
$$
where $(j,l)\in \Nb^{k}$ means that $j\in
\Nb^{k_{1}}$ and $l\in \Nb^{k_{2}}$ with $k_{1}+k_{2}=k$.
Let us define
\be \label{bc1}b_{jl}=i\Omega(j,l)^{-1}a_{jl},\quad c_{jl}=0\quad  
\mbox{when }
\{j_{1},\ldots,j_{k_{1}}\}\neq \{l_{1},\ldots,l_{k_{2}}\}\ee
and
\be\label{bc2}c_{jl}=a_{jl},\quad b_{jl}=0\quad \mbox{when }
\{j_{1},\ldots,j_{k_{1}}\}= \{l_{1},\ldots,l_{k_{2}}\}.\ee
As $\omega$ is strongly non resonant, there exist $\gamma$ and
$\alpha$ such that
$$
|\Omega(j,l)|\geq \gamma\frac{1+S(j,l)}{\mu(j,l)^\alpha}
$$
for all $(j,l)\in \Nb^{k}$ with
$\{j_{1},\ldots,j_{k_{1}}\}\neq \{l_{1},\ldots,l_{k_{2}}\}$.
Thus,
in view of Definitions \ref{po} and \ref{po+}, the polynomial
$$
\chi=\sum_{(j,l)\in \Nb^{k}}b_{j,l}\xi^{(j)}\eta^{(l)},
$$
belongs to $\T_{k}^{\nu',\beta,+}$ while the polynomial
$$
Z=\sum_{(j,l)\in \Nb^{k}}c_{j,l}\xi^{(j)}\eta^{(l)}
$$
belongs to $\T_{k}^{\nu',\beta}$ with $\nu'=\nu +\alpha$. Notice that in this  
non resonant case, \eqref{normal} implies that $Z$ depends only on the actions and thus is in normal form. Furthermore by
construction they
satisfy \eqref{homolo} and \eqref{normal}. Note that the reality of $Q$  is equivalent to the symmetry relation: $\bar a_{jl}=a_{lj}$. Taking
into account that $\Omega_{lj}=-\Omega_{jl}$,
this symmetry remains satisfied for the polynomials
$\chi$ and $Z$. Finally, $\chi$ and $Z$ belong to $\Hs$, since they are homogeneous polynomials (they are their own Taylor expansions) and as a consequence of Proposition \ref{tame} (iii) and (iv) respectively.
\qed  

\medskip
\noindent We can now state the main result of this section:

\begin{theorem}\label{thm:birk}
     Assume that $P$ is a real Hamiltonian belonging to $ \Hs$ for all $s$ large enough and  to  the class $\T^{\nu,\beta}$ for some  
$\nu\geq 0$ and $\beta>0$. Assume
     that $\omega$ is strongly non resonant (cf. Definition \ref{SNR})  
and satisfies \eqref{om} for
     some $\bar d\geq 0$. Then
     for any $r\geq 3$ there exists $s_{0}$ and  for any $s\geq
     s_{0}$ there exists $\U_{s}$, $\V_{s}$ neighborhoods of the origin  
in
     $\Ps$ and
     $\tau_{s} : \V_{s} \to \U_{s}$ a real analytic canonical  
transformation
     which is the restriction to $\V_{s}$ of $\tau := \tau_{s_{0}}$ and
      which puts $H=H_{0}+P$ in normal form up
     to order $r$ i.e.
     $$H\circ \tau= H_{0}+Z+R$$
     with
     \bi
     \item[(a)] $Z$ is a real continuous polynomial of degree $r$ with a
     regular vector field (i.e. $Z\in \mathcal H^s$)  
which only
     depends on the actions: $Z= Z(I)$.
     \item[(b)] $R\in  
\Hs$ is real and $\norma{X_{R}(z)}_{s}\leq
     C_{s}\norma{z}_{s}^{r}$ for all $z\in  
\V_{s}$.
     \item[(c)] $\tau$ is close to the identity: $\norma{\tau
     (z)-z}_{s}\leq C_{s}\norma{z}_{s}^2$ for all $z\in \V_{s}$.
     \ei
\end{theorem}
\proof
The proof is close to the proof of Birkhoff normal form theorem  
stated in \cite{Gre07} or \cite{Bam07}. The main difference has been already pointed out : we have here to check the $\Hs$ regularity  of the Hamiltonian functions at each step, independently of the fact that they belong to $\T^{\nu,\beta}$ (here $P\in \T^{\nu,\beta}$ does not imply $P\in \Hs$).\\
Having fixed some
$r\geq 3$, the idea is to construct iteratively for $ k=
3,\ldots,r$,  a neighborhood $\V_{k}$ of $0$ in $\Ps$ ($s$ large
enough depending on $r$),
a canonical transformation $\tau_k$, defined on $\V_{k}$, an
increasing sequence $(\nu_{k})_{k=3,\ldots,r}$ of positive numbers
and real Hamiltonians $Z_{k}, P_{k+1}, Q_{k+2}, R_{k}$ such that
\begin{equation}
\label{bb.1}
H_{k}:=H\circ\tau_k=H_{0}+Z_{k}+P_{k+1}+Q_{k+2}+R_{k} \,,
\end{equation}
satisfying  the following properties
\bi
\item[(i)] $Z_{k}$ is a polynomial of degree $k$ in
$\T^{\nu_{k},\beta}\cap \Hs$ having a zero of
order 3 (at least) at the origin and $Z_{k}$ depends only on the (new) actions:
$\{Z_{k},I_{j}\}=0$ for all $j\geq 1$.
\item[(ii)] $P_{k+1}$ is a homogeneous polynomial of degree $k+1$ in
$\T_{k+1}^{\nu_{k},\beta}\cap \Hs$.
\item[(iii)] $Q_{k+2}$ is a polynomial of degree $r+1$ in
$\T^{\nu_{k},\beta}\cap \Hs$ having a zero of
order $k+2$ at the origin.
\item[(iv)] $R_{k}$ is a regular Hamiltonian belonging to
$\Hs$ and  having a zero of
order $r+2$ at the origin.
\ei
First we fix $s> \nu_{r}+3/2$ to be sure to be able to apply
Proposition \ref{tame} at each step ($\nu_{r}$ will be defined later
on independently of $s$).
Then we notice that \eqref{bb.1} at order $r$ proves Theorem  
\ref{thm:birk} with
$Z=Z_{r}$ and $R=P_{r+1}+R_{r}$ (since $Q_{r+2}=0$). 
Actually, since $R=P_{r+1}+R_{r}$
belongs to $\mathcal H^{s}$ and has a zero of order $r+1$ at the origin, we can apply Lemma \ref{first} to obtain
\be \label{estimR}
\norma{X_{R}(z)}_{s}\leq C_{s}\norma{z }^r_{s}.
\ee
on $\V\subset \V_{r}$ a neighborhood of $0$ in $\Ps$. \\
At the initial step (which for convenience we will denote the  $k=2$ step), the Hamiltonian $H=H_{0}+P$ has the desired form
(\ref{bb.1}) with $\tau_2=I$, $\nu_{2}=\nu$, $Z_{2}=0$,
$P_{3}$ being the Taylor polynomial of $P$ of degree
$3$, $Q_{4}$ being the Taylor polynomial of $P$ of degree $r+1$
minus $P_{3}$ and
$R_{2}=P-P_{3}-Q_{4}$. We
show now how to go from step $k$ to step $k+1$.\\
We look for $\tau_{k+1}$ of the form $\tau_{k}\circ \phi_{k+1}$,
$\phi_{k+1}$ being the Lie transform associated to a  
homogeneous polynomial $\chi_{k+1}$ of degree $k+1$. 

\medskip 

We decompose $H_{k}\circ \phi_{k+1}$ as follows
\begin{eqnarray}
\label{bb.45} H_{k}\circ \phi_{k+1}&=& H_{0}+Z_{k}+\{H_{0},\chi_{k+1}\}+P_{k+1}
\\
\label{bb.5}
&+& H_{0}\circ \phi_{k+1}-H_{0}-\{H_{0},\chi_{k+1}\}\
\\
\label{bb.6}
&+ & Z_{k}\circ \phi_{k+1} -Z_{k}
\\
\label{bb.7}
&+& P_{k+1}\circ \phi_{k+1}-P_{k+1}
\\
\label{bb.8}
&+&Q_{k+2}\circ \phi_{k+1}
\\
\label{bb.9}
&+& R_{k}\circ \phi_{k+1}\ .
\end{eqnarray}
Using Lemma \ref{homo} above, we choose $\chi_{k+1}$  in $\T_{k+1}^{
\nu'_k,\beta,+}$ in such a way that
\be
\label{bb.10}
\hat{Z}_{k+1}:=\{H_{0},\chi_{k+1}\}+P_{k+1}\ee
is a homogeneous real polynomial of degree $k+1$ in  $\T_{k+1}^{\nu'_k,\beta}$. We put then $Z_{k+1}= Z_k+ \hat{Z}_{k+1}$, which obviously has degree $k+1$ and a zero of order 3 (at least) at the origin, and the right hand side of line \eqref{bb.45} becomes $H_0+Z_{k+1}$. We just recall that $\nu'_k = \nu_k+\alpha$, where $\alpha$ is determined by $\omega$, independently of $r$ and $s$. By Proposition \ref{compo},  the Lie transform associated  
to $\chi_{k+1}$ is well defined
and smooth on
a neighborhood $\V_{k+1}\subset \V_{k}$ and, for $z\in  
\V_{k+1}$ satisfies 
$$
\norma{\phi_{k+1}(z)-z}_s\leq C \norma{z}_s^2.$$
Then from Proposition \ref{bra},  Proposition \ref{compo} and formula  
\eqref{Pphi}, we find that \eqref{bb.6}, \eqref{bb.7}, \eqref{bb.8}
and \eqref{bb.9} are  regular
Hamiltonians
having  zeros of order $k+2$ at the origin. For instance concerning  
\eqref{bb.6}, one has by Taylor formula for any $z\in \V_{k+1}$
$$
Z_{k}\circ \phi_{k+1}(z)  
-Z_{k}(z)=\{Z_k,\chi_{k+1}\}(z)+\int_0^1(1- 
t)\{\{Z_k,\chi_{k+1}\},\chi_{k+1}\}(\Phi^t_{\chi_{k+1}}(z))\ dt
$$
and $\{Z_k,\chi_{k+1}\}$ is a polynomial having a zero of order  
$3+\mbox{degree}(\chi_{k+1})-2=k+2$ while the integral term is a  
regular Hamiltonian having a zero of order $2k+1$. Thus if $2k+1\geq r+2$  
this last term contributes to $R_{k+1}$ and if not, we have to use a  
Taylor formula at a higher order.\\
  Therefore the sum of \eqref{bb.6}, \eqref{bb.7}, \eqref{bb.8}
and \eqref{bb.9} decomposes in $\tilde P_{k+2}+\tilde Q_{k+3}+\tilde  
R_{k+1}$ with
$\tilde P_{k+2}$, $\tilde Q_{k+3}$
and $\tilde R_{k+1}$ satisfying respectively the properties (ii), (iii)  
and (iv) at rank $k+1$ (with
$\nu_{k+1}=k\nu'_{k}+ \nu_{k}+k+2$).\\
Concerning the term \eqref{bb.5}, one has to proceed differently since  
$H_0$ does not belong to the $\Hs$. \\
First notice that by  the homological equation \eqref{bb.10} one has  $\{H_{0},\chi_{k+1}\}=Z_{k+1}-Z_{k}-P_{k+1}$. By construction $Z_k$ and $P_{k+1}$ belong to $\Hs$. On the other hand, by Lemma \ref{homo}, $Z_{k+1}\in \T^{\nu'_k,\beta}_{k+1}$ and is in normal form (i.e. it depends only on the action variables). Thus by Proposition \ref{tame} assertion (iv), one concludes that $Z_{k+1}\in\Hs$. Therefore we have proved that   $\{H_{0},\chi_{k+1}\}\in \Hs$. \\
Now we use the Taylor formula at order one to get
$$
H_{0}\circ  
\phi_{k+1}(z)-H_{0}(z)=\int_0^1\{H_0,\chi_{k+1}\}(\Phi^t_{\chi_{k+1}}(z))\ dt.
$$
But we know from the proof of Proposition \ref{compo} that  $\Phi^t_{\chi_{k+1}}: \V_{k+1}\to \Ps$ for all $t\in[0,1]$.  Therefore  $H_{0}\circ  
\phi_{k+1}-H_{0}\ \in \ \Hs$ and thus \eqref{bb.5} defines a regular Hamiltonian.\\
Finally we use again the Taylor formula and the homological equation to write
\begin{equation*}\begin{split}
H_{0}\circ  
\phi_{k+1}(z)-H_{0}(z)-&\{H_{0},\chi_{k+1}\}(z)=\\
&\int_0^1(1-t)\{Z_{k+1}- 
Z_k-P_{k+1},\chi_{k+1}\}(\Phi^t_{\chi_{k+1}}(z))\ dt
\end{split}\end{equation*}
and , since $Z_{k+1}-Z_k-P_{k+1}$ belongs to $\T_{k+1}^{
\nu'_k,\beta}$ and $\chi_{k+1}\in \T_{k+1}^{
\nu'_k,\beta,+}$  we conclude by Proposition \ref{bra} that  $H_{0}\circ  
\phi_{k+1}-H_{0}-\{H_{0},\chi_{k+1}\}\in \T^{
\nu_{k+1},\beta}$. Finally we use Proposition \ref{compo} 
to decompose it in $\hat P_{k+2}+\hat Q_{k+3}+\hat R_{k+1}$ with
$\hat P_{k+2}$, $\hat Q_{k+3}$
and $\hat R_{k+1}$ satisfying respectively the properties (ii), (iii)  
and (iv) at rank $k+1$. The proof is achieved defining $P_{k+2}=\hat  
P_{k+2}+\tilde P_{k+2}$, $ Q_{k+3}=\hat Q_{k+3}+\tilde Q_{k+3}$ and $  
R_{k+1}=\hat R_{k+1}+\tilde R_{k+1}$.
  \endproof

\section{Dynamical consequences}\label{dyn}
\subsection{Nonlinear harmonic oscillator in one dimension}\label{nls1d}
We recall the notations of the introduction.
The quantum harmonic oscillator $T=-\frac{d^2}{dx^2}+x^2$ is diagonalized in the Hermite basis $(\phi_j)_{j\in\Nb}$:
\begin{align*}T\phi_j &=(2j-1)\phi_j, \quad j\in \Nb\\
\phi_{n+1} &=\frac{H_n(x)}{\sqrt{2^n n!}}e^{-x^2/2} , \quad n\in \N
\end{align*}
where $H_n(x)$ is the $n^{th}$ Hermite polynomial relative to the weight $e^{-x^2}$ :
$$
\int_\R e^{-x^2} H_m(x)H_n(x)dx=2^n n! \sqrt\pi \delta_{nm}.$$
In this basis, the Hermite multiplier is given by
\be \label{M1}
M\phi_j=m_j \phi_j
\ee
where $(m_j)_{j\in \Nb}$ is a bounded sequence of real number.
For any $k\geq 1$, we define the class
\be \label{Wk}
\W_k=\{   (m_j)_{j\in \Nb}\mid \mbox{ for each }j,\ m_j=\frac{\tilde m_j}{j^k}\mbox{ with } \tilde m_j\in [-1/2,1/2]\}
\ee
that we endow with the product probability measure. In this context the frequencies, i.e. the eigenvalues of $T+M= -d^2/dx^2 +x^2+M$  are given by
$$\omega_j= 2j-1+m_j=2j-1 +\frac{\tilde m_j}{j^k}, \quad j\in \Nb.$$

\begin{proposition}\label{res}
There exists a set $F_{k} \subset \W_{k}$
whose measure equals $1$ such that if $m=(m_j)_{j\in \Nb}\in F_{k}$ then the 
frequency vector
$(\omega_{j})_{j\geq 1}$ is strongly non-resonant (cf. Definition \ref{SNR}).
\end{proposition}
\proof
First remark that it suffices to prove that the 
frequency vector
$(\omega_{j})_{j\geq 1}$ is  non resonant in the sense of Definition \ref{NR}. Actually, if we prove that \eqref{A.1} is satisfied for given constants $\delta'$ and $\gamma'$ then 
  if $S(j)<r\mu(j)$ 
 $$    \left|\omega_{j_1}+\cdots+\omega_{j_{i}}-\omega_{j_{i+1}}-\cdots 
    -\omega_{j_{r}} \right|\geq \frac{\gamma'}{\mu(j)^{\delta'}}
  \geq \frac{\gamma'}{r+1}\frac{1+S(j)}{\mu(j)^{\delta'+1}}
$$
and thus \eqref{A.2} is satisfied with $\delta =\delta' +1$ and $\gamma= \frac{\gamma'}{r+1}$. Now if $S(j)\geq r\mu(j)$ then  use
\begin{equation} 
    \label{A.3} 
    \left|\omega_{j_1}+\cdots+\omega_{j_{i}}-\omega_{j_{i+1}}-\cdots 
    -\omega_{j_{r}} \right|\geq S(j)-(r-2)\mu(j),  
    \end{equation} 
to conclude that
$$
     \left|\omega_{j_1}+\cdots+\omega_{j_{i}}-\omega_{j_{i+1}}-\cdots 
    -\omega_{j_{r}} \right| \geq \frac{2}{r}S(j)
   \geq \frac{\gamma'}{r+1}\quad \frac{1+S(j)}{\mu(j)^{\delta'+1}}
    $$
    provided $\gamma'$ is small enough.\\
The proof that there exists a set $F_{k} \subset \W_{k}$
whose measure equals $1$ such that if $m=(m_j)_{j\in \Nb}\in F_{k}$ then the 
frequency vector
$(\omega_{j})_{j\geq 1}$ is  non resonant is exactly the same as the proof of Theorem 5.7 in \cite{Gre07}. So we do not repeat it here (see also \cite{BG04}).
\endproof

In equation \eqref{1} with $d=1$, the Hamiltonian perturbation reads
\be \label{P}
P(\xi,\eta)=\int_{\R}g(\xi(x),\eta(x))dx\ee
where $g\in C^\infty ( \C^2,\C)$, $\xi(x)=\sum_{j\geq 1}\xi_{j}\phi_{j}(x)$,
$\eta(x)=\sum_{j\geq 1}\eta_{j}\phi_{j}(x)$ and\\ $((\xi_{j})_{j\geq 1}, 
(\eta_{j})_{j\geq 1})\in \Ps$.  We first  check that $P$ belongs to $\Hs$ for $s$ large enough. 

\begin{lemma}\label{PHs}
Let $P$ given by \eqref{P} with $g\in C^\infty (\U, \C)$, $\U$ being a neighborhood of $0$ in $\C^2$, $g$ real i.e. $g(z,\bar z)\in \R$ and $g$ having a zero of order at least 3 at the origin. Then $P\in \Hs$ for all $s>1/2$.
\end{lemma}
\proof
One computes
$$
\frac{\partial P}{\partial \xi_j}(\xi,\eta)=\int_{\R}\partial_1 g(\xi(x),\eta(x)) \phi_j(x)dx
$$
and
$$
\frac{\partial P}{\partial \eta_j}(\xi,\eta)=\int_{\R}\partial_2 g(\xi(x),\eta(x)) \phi_j(x)dx\,.
$$
In the same way, we have
\begin{eqnarray}
\lefteqn{\frac{\partial^{l+r} P}{\partial \xi_{j_1}\dots\partial \xi_{j_l}\partial \eta_{k_1}\dots\partial \eta_{k_r}}(\xi,\eta)}&&\nonumber\\
&=& \int_\R \partial_1^l \partial_2^r g (\xi(x),\eta(x))\phi_{j_1}(x)\dots\phi_{j_l}(x)\phi_{k_1}(x)\dots\phi_{k_r}(x)dx\,.\label{dlrp}
\end{eqnarray}
Since $g$ is a $C^\infty$ function, all these partial derivatives are continuous from $\Ps$ to $\C$, and the corresponding differentials $(\xi,\eta) \rightarrow D^{l+r}P(\xi,\eta)$ are continuous from $\Ps$ to the space of $l+r$-linear forms on $\Ps$. We get moreover
\begin{equation*} \begin{split}
\norma{X_P(\xi,\eta)}_s^2&= \sum_{j\geq 1}\vaj^{2s} \Va{\int_{\R}\partial_1 g(\xi(x),\eta(x)) \phi_j(x)dx}^2\\
&+\sum_{j\geq 1}\vaj^{2s} \Va{\int_{\R}\partial_2 g(\xi(x),\eta(x)) \phi_j(x)dx}^2.
\end{split}
\end{equation*}
Therefore, to check that $z\mapsto X_P(z)$ is a regular function from a neighborhood of the origin in $\Ps$ into $\Ps$, it suffices to check that the functions $x\mapsto \partial_1 g(\xi(x),\eta(x))$ and $x\mapsto \partial_2 g(\xi(x),\eta(x))$ are in $\tilde H^s$ provided $\xi(x)$ and $\eta(x)$ are in $\tilde H^s$. So it remains to prove that functions of the type \\$x\mapsto \va{x}^i\partial_1^{l+1}\partial_2^mg(\xi(x),\eta(x))(\xi^{(l_1)}(x))^{\alpha_1}\cdot (\xi^{(l_{k_1})}(x))^{\alpha_{k_1}}(\eta^{(m_1)}(x))^{\beta_1}\cdots (\eta^{(m_{k_2})}(x))^{\beta_{k_2}}$ are in $ L^2(\R)$ for all $0\leq i+l+m\leq s$, $0\leq i+l_j\leq s$, $0\leq i+m_j\leq s$.
But this is true because 
\begin{itemize}
\item $g$ is a $C^\infty$ function, $\xi$ and $\eta$ are bounded functions and thus 
$x\mapsto \partial_1^{l+1}\partial_2^mg(\xi(x),\eta(x))$ is bounded
\item 
$\tilde H^s$ is an algebra for $s>1/2$ and thus $x\mapsto \va{x}^k\xi^{(l)}(x)\eta^{(m)}(x)\ \in \ L^2(\R)$ for all $0\leq k+l+m\leq s$.
\item $ |\partial_1 g( \xi(x),\eta(x))|,\  |\partial_2 g( \xi(x),\eta(x))|\leq C (|\xi(x)|+|\eta(x)|)^2$
for some uniform constant $C>0$ and thus $x\mapsto \va{x}^k\partial_1g(\xi(x),\eta(x))\ \in \ L^2(\R)$ for all $0\leq k\leq s$.
\end{itemize}
There remains to prove the same properties concerning the Taylor homogeneous polynomial $P_m$ of $P$ at any order $m$, computed at $(0,0)$. From (\ref{dlrp}), we get
\begin{equation*}
P_m = \frac{1}{m!}\int_\R \sum_{l+r=m} \partial_1^l \partial_2^r g (0,0)\sum_{j,k} \xi_{j_1}\phi_{j_1}(x)\dots\xi_{j_l}\phi_{j_l}(x)\eta_{k_1}\phi_{k_1}(x)\dots\eta_{k_r}\phi_{k_r}(x)dx\,,
\end{equation*}
hence $P_m$ can be computed directly from formula (\ref{P}), replacing $g$ by its Taylor homogeneous polynomial $g_m$ of order $m$ :
\begin{equation*}
g_m(\xi(x),\eta(x))  = \frac{1}{m!} \sum_{l+r=m} \partial_1^l \partial_2^r g (0,0)\xi(x)^l\eta(x)^r ,,
\end{equation*}
and this gives the statement, since $g_m$ satisfies the same properties as $g$. 
\endproof

The fact that $P$ belongs to the class $\T^{\nu,\beta}$ is directly related to 
the distribution of the $\phi_{j}$'s. Actually we have
\begin{proposition}\label{prop:phi} Let $\nu>1/8$ and $0 \leq \beta \leq \frac{1}{24}$.
For each $k\geq 1$ and for each $N\geq 0$ there exists 
$c_N>0$ such that for all $j\in \Nb^k$
\be \label{phi}
\Va{\int_{\R}\phi_{j_{1}}\ldots\phi_{j_{k}} dx}\leq c_N 
\frac{\mu(j)^{\nu}}{C(j)^\beta}A(j)^N.\ee
As a consequence,  any $P$ of the general form \eqref{P}  is in the class $\T^\nu$. 
\end{proposition}
The proof will be done in the multidimensional case in the next section (cf. Proposition \ref{prop:phid}).

\medskip

We can now apply our Theorem \ref{thm:birk} to obtain 
\begin{theorem}\label{thm:dyn} 
 Assume that  $M\in F_{m}$ defined in Proposition \ref{res} and that $g\in C^\infty ( \C^2,\C)$ is real i.e. $g(z,\bar z)\in \R$ and has a zero of order at least 3 at the origin. For any $r\geq 3$ there exists $s_{0}(r)$ an integer such that for any  
 $s\geq s_{0}(r)$, there exist 
 $\e_{0}>0$ and $C>0$ such that if $\norma{\psi_0}_{\tilde H^{2s}}=\e <\e_{0}$ the equation
\be \label{eq}i\psi_t=(-\Delta +x^{2} +M)\psi +\partial_2 g(\psi,\bar \psi), \quad x\in \R,\  t\in \R \ee
 with Cauchy data $\psi_0$
has a unique solution 
$ \psi\in 
C^1((-T_{\epsilon},T_{\epsilon}),\tilde H^{2s})$ with
\be \label{t}
T_{\epsilon}\geq C\epsilon^{-r} .
\ee
Moreover 
\begin{itemize}
\item[(i)] $\Vert{\psi(t,\cdot )}\Vert_{\tilde H^{2s}}\leq 2\epsilon$ for any $t\in (-T_{\epsilon},T_{\epsilon})$. 
\item[(ii)] 
$\sum_{j\geq 1}j^{2s}|\va{\xi_j(t)}^2-\va{\xi_j(0)}^2|\leq {\e^{3}} $ for any $t\in (-T_{\epsilon},T_{\epsilon})$\\
where $\va{\xi_j(t)}^2$, $j\geq 1$ are the actions of $\psi(t,\cdot)=\sum \xi_j(t)\phi_j$ .
 \item[(iii)] there exists a torus 
$\mathcal T_{0}\subset \tilde H^{2s}$ such that, 
$$\mbox{dist}_{2s}(\psi(t,\cdot),\mathcal T_{0})\leq C\e^{r_1/2} \quad \mbox{ for } \vat\leq \epsilon^{-r_2} $$
where $r_1+r_2=r+3$ and  $\mbox{dist}_{2s}$ denotes the distance on $\tilde H^{2s}$ associated with 
the norm $\norma{\cdot}_{\tilde H^{2s}}$.
\end{itemize}
\end{theorem}
\proof
Let $\psi_0 =\sum \xi_j(0) \phi_j$ and denote $z_0= (\xi(0),\bar\xi(0))$. Notice that if $\psi_0\in \tilde H^{2s}$ with  $\norma{\psi_0}_{\tilde H^{2s}}=\epsilon$ then $z_0 \in \Ps$ and $\norma{z_0}_s=\epsilon$. Denote by $z(t)$ the solution of the Cauchy problem $\dot z= X_H(z)$, $z(0)=z_0$, where  $H=H_0+P$ is the Hamiltonian function associated to the equation \eqref{eq} written in the Hermite decomposition $\psi(t) =\sum \xi_j(t) \phi_j$, $z(t)=(\xi(t), \bar \xi(t))$. We note that, since $P$ is real, $z$ remains a real point of $\Ps$ for all $t$ and that $\norma{\psi(t)}_{\tilde H^{2s}}=\norma{z(t)}s$.\\
Then we denote by $z'=\tau^{-1}(z)$ where $\tau :\V_{s} \to \U_{s}$ is the 
transformation given by Theorem \ref{thm:birk} (so that $z'$ denotes the normalized coordinates) associated to the order $r+2$ and $s\geq s_0(r+2)$ given by the same Theorem. We note that, since the transformation $\tau$ is generated by a real Hamiltonian, $z'(t)$ is still a real point.\\
 Let $\e_0>0$ be such that $B_{2\e_0}\subset  \V_{s}$ and take $0<\e<\e_0$. We assume that $\norma{z(0)}_s=\norma{\psi_0}_{\tilde H^{s}}=\e$. For $z=(\xi,\eta)\in \Ps$ we define
$$N(z):=2\sum_{j=1}^\infty j^{2s} I_{j}(\xi,\eta)$$
 where we recall that $I_j(\xi,\eta)=\xi_j\eta_j$.
We notice that for a real point $z=(\xi,\bar\xi)\in \Ps$, 
$$N(z)=\norma{z}_{s}^2.$$
Thus in particular we have\footnote{That is precisely at this point  that we need to work with real Hamiltonians. The Birkhoff normal form theorem is essentially algebraic and does hold for complex Hamiltonians.  }
$$ N(z(t))=\norma{z(t)}_{s}^2\quad \mbox{ and }N(z'(t))=\norma{z'(t)}_{s}^2.$$
Using that $Z$ depends
only on the normalized actions, we have
\be \label{ca}
\dot N(z')=\{N ,H\circ \tau \}\circ \tau^{-1}(z)= \{N,
R\}(z').\ee
Therefore as far as $\norma{z(t)}_s\leq 2\e$, and thus $z(t)\in \V_{s}$,  by assertion (c) of Theorem \ref{thm:birk}, $\norma{z'(t)}_s\leq C\e$
and using   \eqref{ca} and assertion (b) of Theorem \ref{thm:birk} (at order $r+2$) we get 
$$ \Va{N(z'(t))-N(z'(0))}\leq  \Va{\int_0^t\{N,
R\}(z'(t'))dt'}\leq Ct\norma{z'(t)}_s^{r+3}\leq Ct \e^{r+3 }.$$
In particular, as far as $\norma{z(t)}_s\leq 2\e$ and $\vat\leq C\e^{-r}$
$$ \Va{N(z'(t))-N(z'(0))}\leq C \e^3.$$
Therefore using again assertion (c) of Theorem \ref{thm:birk}, we obtain
$$ \Va{N(z(t))-N(z(0))}  \leq C\e^3$$
which, choosing $\e_0$ small enough, leads to $\norma{z(t)}_s\leq 3/2\ \e$ as long as $\norma{z(t)}_s\leq 2\e$ and $\vat\leq C\e^{-r}$. Thus   \eqref{t} and assertions (i) follow  by a continuity argument.

\medskip

To prove assertion (ii) we recall the notation
$I_j(z)=I_j(\xi,\eta)=\xi_j\eta_j$ for the actions associated to $z=(\xi,\eta)$. 
Using that $Z$ depends 
only on the actions, we have
$$
\{I_j\circ \tau^{-1},H\}(z)=\{I_j ,H\circ \tau \}\circ \tau^{-1}(z)= \{I_j, 
R\}(z').$$ 
Therefore, we get in the normalized 
coordinates
$$
\frac d {dt}I_j(\xi',\eta')=-i\xi'_j\frac{\partial R}{\partial
  \eta_j}+i\eta'_j\frac{\partial R}{\partial \xi_j}\  
$$
and thus 
\begin{eqnarray*}
\sum_jj^{2s}\left| \frac d {dt} I_j(\xi',\bar\xi')\right|&=&
  \sum_jj^{2s}\left|-\xi'_j\frac{\partial R}{\partial
  \eta_j}+\bar\xi'_j\frac{\partial R}{\partial \xi_j} \right| 
\\
&\leq&
  \left(\sum_{j}j^{2s}(|\xi'_j|^2+|\xi'_j|^2) \right)^{1/2}
  \left(\sum_{j}j^{2s}\left(\left|\frac{\partial R}{\partial
  \eta_j}\right|^2+\left|\frac{\partial R}{\partial \xi_j}\right|^2\right)
  \right)^{1/2}
\end{eqnarray*}
which leads to 
\begin{equation}\label{n.ks}
    \sum_jj^{2s}\left| \frac d {dt} I_j(z')\right|\leq 
    \norma{z'}_{s}\norma{X_{R}(z')}_{s} \leq  \norma{z'}_{s}^{r+3}.
    \end{equation}
Thus, recalling that $I_j(\xi',\bar \xi')=\va{\xi'_j}^2$   we get
\be \label{Inew}\sum_{j\geq 1}j^{2s}\Va{\va{\xi'_j(t)}^2-\va{\xi'_j(0)}^2}\leq {\e^{3}} \mbox{ for any }\vat\leq C \e^{-r}.\ee
On the other hand, using (i) and assertion (c) of Theorem \ref{thm:birk}, for any $\vat \leq C \e^{-r}$, one has
$$
\sum_{j\geq 1}j^{2s}\Va{\va{\xi_j(t)}^2)-\va{\xi'_j(t)}^2}\leq \sum_{j\geq 1}j^{2s}(\va{\xi_j(t)}+\va{\xi'_j(t)})\Va{\xi_j(t)-\xi'_j(t)}\leq C \e^3.$$
Combining this last relation with \eqref{Inew}, assertion (ii) follows.

\medskip

To prove (iii), let $\bar I_{j}= I'_{j}(0)$ be the initial actions in the normalized 
coordinates and define the smooth torus
$$
\Pi _0:=\left\{z\in\mathcal P_{s}\ :\ I_j(z)=\bar I_j\ ,j\geq 1 \right\}
$$
and its image in $\tilde H^{s}$
$$
\mathcal T_0= \{u\in \tilde H^{s}\ :\ u=\sum \xi_j \phi_j \mbox{ with } \tau(\xi,\bar \xi)\in \Pi _0\}.
$$
We have
\begin{equation}
\label{de2}
d_{s}(z(t),\mathcal T_0)\leq \left[\sum_{j}j^{2s}\left| \sqrt{I'_j(t)}-\sqrt{\bar
  I_j} \right|^2 \right]^{1/2}
\end{equation}
where $d_s$ denotes the distance in $\Ps$ associated to $\norma{\cdot}_s$.\\
Notice that for $a,b\geq 0$,
$$
\left|\sqrt{a}-\sqrt{b}\right|\leq \sqrt{\va{a-b}}\ .
$$ 
Thus, using
(\ref{n.ks}), we get
\begin{eqnarray*}
\left[d_{s}(z(t),\mathcal T_0) \right]^2 &\leq &\sum_{j}j^{2s}|I'_j(t)-
  I'_j(0)|\\ &\leq &|t| \sum_{j}j^{2s}|\dot I'_j(t)| \\
  &\leq & \frac 
  1{\epsilon^{r_1}}\norma{z'}_s\norma{X_{R}(z')}_{s}\\
  &\leq & C\frac 
  1{\epsilon^{r_1}}\epsilon^{r+3}\leq C \epsilon^{r+3-r_1}.
  \end{eqnarray*}
which gives (ii).

\endproof

\subsection{Multidimensional nonlinear harmonic oscillator}\label{nlsdd}
\subsubsection{Model}
The spectrum of the d-dimensional harmonic oscillator $$T=-\Delta + \va{x}^2=-\Delta +x_1^2+\cdots+x_d^2$$ is  the sum of $d$-copies of the odd integers set, i.e. the spectrum of $T$ equals $\N_d$ with
\be\label{Nd}
\N_d= \left\{ \begin{array}{c}2\N \setminus \{0,2,\cdots,d-2\}\mbox{ if d is even}\\
2\N+1\setminus \{1,3,\cdots,d-2\} \mbox{ if d is odd}.
\end{array}\right.\ee
For $j\in \Nd$ we denote the associated eigenspace $E_j$ which dimension is
$$d_j=\sharp\{ (i_1, \cdots ,i_d)\in(2\N+1)^d \mid i_1+\cdots+i_d=j \}.$$
We denote $\{\Phi_{j,l}$, $l=1,\cdots,d_j\}$, the basis of $E_j$ obtained by  $d$-tensor product  of Hermite functions:   $\Phi_{j,l}=\phi_{i_1} \otimes\cdots \phi_{i_d}$ with $i_1+\cdots+i_d=j $ .

The Hermite multiplier $M$ is defined on the basis $(\Phi_{j,l})_{j\in \Nd,l=1,\cdots,d_j}$ of $L^2(\R^d)$ by
\be\label{M}
M \Phi_{j,l}=m_{j,l}\Phi_{j,l}
\ee
where $(m_{j,l})_{j\in \Nd,l=1,\cdots,d_j}$ is a bounded sequence of real numbers.

The linear part of \eqref{1} reads
$$H_0= -\Delta +x^2+M .$$
$H_0$ is still diagonalized by $(\Phi_{j,l})_{j\in \Nd,l=1,\cdots,d_j}$ and  the spectrum of $H_0$ is
\be \label{spectre}
\sigma (H_0) = \{j+m_{j,l} \mid {j\in \Nd,l=1,\cdots,d_j}   \}
\ee
For simplicity, we will focus on the case $m_{j,l}=m_j$ for all $l=1,\cdots,d_j$.
In this case we have $\sigma (H_0) = \{j+m_{j} \mid {j\in \Nd}   \}$ and, as a consequence of Proposition \ref{res},
\begin{proposition}\label{resd}
There exists a set $F_{k} \subset \W_{k}$
whose measure equals $1$ such that if $m=(m_j)_{j\in \Nb}\in F_{k}$ then the 
frequency vector
$(\omega_{j,i})_{j\in \Nd,i=1,\cdots,d_j}$ satisfies the following:\\
for  any $r\in\Nb$, there are $ \gamma >0$ and 
    $\delta >0$ such that for any $j\in \Nd^{r}$, any $l\in \{1,\cdots,d_{j_1}\}\times \cdots \times
 \{1,\cdots,d_{j_r}\}$     and any $1\leq i\leq 
    r$, one has 
    \begin{equation} 
    \label{A.2d} 
    \left|\omega_{j_1,l_1}+\cdots+\omega_{j_{i},l_i}-\omega_{j_{i+1},l_{i+1}}-\cdots 
    -\omega_{j_{r},l_r} \right|\geq \gamma \frac{1+S(j)}{\mu(j)^{\delta}}  
    \end{equation} 
except if $\{j_{1},\ldots,j_{i}\}=\{j_{i+1},\ldots,j_{r}\}$.    
\end{proposition}
Concerning the product of eigenfunctions we have,
\begin{proposition}\label{prop:phid}
Let $\nu>d/8$. For any $k\geq 1$ and any $N\geq 1$ there exists $c_N>0$ such that for any $j\in \Nd^{k}$, any $l\in \{1,\cdots,d_{j_1}\}\times \cdots \times
 \{1,\cdots,d_{j_k}\}$
\be \label{phid}
\Va{\int_{\R^d}\Phi_{j_{1},l_1}\ldots\Phi_{j_{k},l_k} dx}\leq c_N 
\frac{\mu(j)^{\nu}}{C(j)^{\frac 1 {24}}}A(j)^N.\ee
\end{proposition}
Notice that this condition does not distinguish between modes having the same energy.
\proof We use the approach developed  in \cite{Bam07} Section 6.2. 
The basic idea lies in the following commutator lemma: Let $A$ be a linear operator which maps $D(T^k)$ into itself and define the sequence of operators
$$
A_N:=[T,A_{N-1}], \quad A_0:=A$$
then (\cite{Bam07} Lemma 7) for any $j_1\neq j_2$ in $\Nd$, any $0\leq l_1\leq d_1$, $0\leq l_2\leq d_2$ and any $N\geq 0$
$$
\va{\langle A\Phi_{j_2,l_2},\Phi_{j_1,l_1}\rangle}\leq \frac{1}{\va{{j_1}-{j_2}}^N}\va{\langle A_N\Phi_{j_2,l_2},\Phi_{j_1,l_1}\rangle}.$$
Let $A$ be the  operator  given by the multiplication by the function $\Phi=\Phi_{j_3,l_3}\cdots\Phi_{j_k,l_k}$ then by an induction argument
$$
A_N=\sum_{0\leq\va{\alpha}\leq N}C_{\alpha,N}D^\alpha$$
where
$$
C_{\alpha,N}= \sum_{0\leq \va{\beta}\leq 2N-\va{\alpha}} V_{\alpha,\beta,N}(x) D^\beta \phi$$
and $V_{\alpha,\beta,N}$ are polynomials of degree less than $2N$. Therefore one gets
\begin{align}\begin{split}\label{3.5}
\Va{\int_{\R^d}\Phi_{j_{1},l_1}\ldots\Phi_{j_{k},l_k} dx}&\leq  \frac{1}{\va{{j_1}-{j_2}}^N}\norma{A_N \Phi_{j_2,l_2}}_{L^2}\\
&\leq C   \frac{1}{|j_{1}-j_{2}|^N } 
\sum_{0\leq \va{\alpha} \leq  N} \sum_{0\leq \va{\beta} \leq 2N-\va{\alpha}} ||V_{\alpha,\beta,N}   D^\beta \phi D^\alpha \Phi_{j_2,l_2} ||_{L^2}  \\
& \leq C  \frac{1}{|j_{1}-j_{2}|^N } \sum_{0\leq \va{\alpha} \leq N} \sum_{0\leq \va{\beta} \leq 2N-\va{\alpha}}  \norma{\Phi_{j_2,l_2}}_{\va{\alpha}}  ||\Phi ||_{{\nu_0+\va{\beta}}}                        
\end{split}\end{align}
where we used in the last estimate (in this proof, $\norma{f}_s=\norma{f}_{H^s(\R^d)}$, the standard  Sobolev norm)
$$\forall \nu_0 > d/2 \qquad ||fg||_{L^2} \leq C_{\nu_0} ||f||_{{\nu_0}} ||g||_{L^2} .            $$
We now estimate $||\Phi ||_{{\nu_0+\va{\beta}}} $.
First notice that, since $T\Phi_{j,l}=j\Phi_{j,l}$, one has for all $s\geq 0$
\be\label{phi-s}
\norma{\Phi_{j,l}}_s\leq Cj^{s/2}.\ee
Then we recall that the Hermite eigenfunctions are uniformly bounded, and in fact (see \cite{Szeg75} or \cite{KoTa05})
\begin{equation}
\label{szego} 
   ||\phi_j||_{L^\infty}\leq C j^{-1/12}   ,
\end{equation}
and thus, since $\Phi_{j,l}=\phi_{i_1} \otimes\cdots \phi_{i_d}$ with $i_1+\cdots+i_d=j$,
we deduce
\begin{equation}
\label{szegod} 
   ||\Phi_{j,l}||_{L^\infty}\leq C_d j^{-1/12}   
\end{equation}
with $C_d= C d^{1/12}$.  
Thus using  tame estimates (see for instance \cite{Tay91})
$$
\norma{uv}_s\leq C(\norma{u}_s\norma{v}_{L^\infty}+\norma{v}_s\norma{u}_{L^\infty})
$$
combined with  \eqref{phi-s} and \eqref{szegod}, we get
for $j_3\geq \cdots\geq j_k$,
\be\label{estimphi}
\norma{\Phi}_s\leq C j_3^{s/2}.
\ee
Inserting \eqref{phi-s} and \eqref{estimphi} in \eqref{3.5} we get
\begin{align*}
\Va{\int_{\R^d}\Phi_{j_{1},l_1}\ldots\Phi_{j_{k},l_k} dx}&\leq C  \frac{1}{|j_{1}-j_{2}|^N } \sum_{0\leq \va{\alpha} \leq N} \sum_{0\leq \va{\beta} \leq 2N-\va{\alpha}}  j_{2}^{\va{\alpha}/2} j_{3}^{(\nu_0+\va{\beta})/2)}      \\
&\leq C  \frac{1}{|j_{1}-j_{2}|^N }   \sum_{0\leq \va{\alpha} \leq N} j_{2}^{\va{\alpha}/2} j_{3}^{\nu_0/2+N-\va{\alpha}/2}\\
&\leq       C\frac{1}{|j_{1}-j_{2}|^N} j_{3}^{N+\nu_0/2} \left( \frac{j_{2} }{j_{3}}\right)^{N/2}\\
 &= C\frac{1}{|j_{1}-j_{2}|^N}j_{3}^{\nu_0/2}{(j_{2}j_{3})^{N/2} }.     
\end{align*}
Now remark that if $\sqrt{j_{2}j_{3}}\leq  |j_{1}-j_{2}|$ then the last estimate implies
\be\label{presque} \Va{\int_{\R^d}\Phi_{j_{1},l_1}\ldots\Phi_{j_{k},l_k} dx}  \leq C j_{3}^{\nu_0/2}  \frac{(j_{2}j_{3})^{N/2}  }{(\sqrt{j_{2}j_{3}}+|j_{1}-j_{2}|)^N }  =C\mu(j)^{\nu/2} A(j)^N    \ee
while if $\sqrt{j_{2}j_{3}} > |j_{1}-j_{2}|$ then $A(j)\geq 1/2$ and thus \eqref{presque} is trivially true.

On the other hand, using \eqref{szegod} one has 
$$
\Va{\int_{\R^d}\Phi_{j_{1},l_1}\ldots\Phi_{j_{k},l_k} dx} \leq C j_{1}^{-1/12}=C \ C(j)^{-1/12} . $$
Combining this estimate with \eqref{presque} one gets for all $N\geq 1$
$$\Va{\int_{\R^d}\Phi_{j_{1},l_1}\ldots\Phi_{j_{k},l_k} dx}\leq c_N 
\frac{\mu(j)^{\nu}}{C(j)^{\frac 1 {24}}}A(j)^N$$
with $\nu=\frac{\nu_0}{4}$.
\endproof

\subsubsection{Result}
We first generalize the normal form theorem to a context adapted to the multidimensional case.
We follow the 
presentation of Section \ref{birk} and only focus on the new 
features.

Let $s\geq 0$, we consider the phase 
space $\Qs=\ls \times \ls$ with
$$\ls=\{(a_{j,l})_{j\in \Nd,\,
1\leq l\leq d_{j}} \mid \sum_{j\in \Nd 
}\vaj^{2s}\sum_{l=1}^{d_{j}}\va{a_{j,l}}^2 <\infty \}
$$
that we endow with the standard norm and the standard symplectic 
structure as for $\Ps$ in Section \ref{model}. Writing $\psi=\sum \xi_{j,l}\Phi_{j,l}$, 
$\bar \psi=\sum \eta_{j,l}\Phi_{j,l}$ with
$(\xi,\eta)\in \Qs$, we note that $\psi\in \tilde H^{2s}$ if and only if $\xi\in \ls$. The linear part of the multidimensional version of the linear part  of \eqref{1} reads
$$
H_{0}(\xi,\eta)=\frac{1}{2}\sum_{j\in \Nd 
}\sum_{l=1}^{d_{j}} \omega_{j,l} \xi_{j,l}\eta_{j,l}.
$$
For $j\geq 1$, we define
$$
J_{j}(\xi,\eta)=\sum_{l=1}^{d_{j}}\xi_{j,l}\eta_{j,l} \ .
$$
Using notations of Section \ref{model}, we define
the class $\Tt_{k}^{\nu}$ of  homogeneous polynomials of 
degree $k$ on $\Qs$
$$Q(\xi, \eta) \equiv Q(z)=\sum_{j\in 
\Nd^k}\sum_{l_{1}=1}^{d_{j_{1}}}\ldots\sum_{l_{m}=1}^{d_{j_{k}}}
a_{j,l}z_{j_{1},l_{1}}\ldots z_{j_{k},l_{k}}$$
such that for each $N\geq 1$, there exists a constant $C>0$ such that for all 
$j,l$
$$\va{a_{j,l}}\leq C \frac{\mu(j)^{\nu}}{C(j)^{1/24}}A(j)^N.$$
Then, following Definition \ref{T} we define a corresponding class 
$\Tt^\nu$ of $C^\infty$  Hamiltonians on $\Qs$ having their Taylor polynomials in $\Tt_{k}^{\nu}$. Similarly, following Definition \ref{defHs}, we also define $\Hsd$ the class of real Hamiltonians $P$ satisfying $P,P_k\in C^\infty(\U_s,\C)$ and $X_P , X_{P_k}\in C^\infty(\U_s,\Qs)$ for some $\U_s\subset \Qs$ a neighborhood of the origin and for all $k\geq 1$ (as before $P_k$ denotes the Taylor polynomial of $P$ of degree $k$).\\
In equation \eqref{1}, the Hamiltonian perturbation reads
\be \label{Pd}
P(\xi,\eta)=\int_{\R^d}g(\xi(x),\eta(x))dx\ee
where $g$ is $ C^\infty $ on a neighborhood of $0$ in $\C^2$, $\xi(x)=\sum_{j\geq 1}\xi_{j}\phi_{j}(x)$,
$\eta(x)=\sum_{j\geq 1}\eta_{j}\phi_{j}(x)$ and $((\xi_{j})_{j\geq 1}, 
(\eta_{j})_{j\geq 1})\in \Ps$.  As in the one dimensional case (cf. Lemma \ref{PHs}),  $P$ belongs to $\Hsd$ for $s$ large enough ($s>d/2$) and using Proposition \ref{prop:phid}, $P$ belongs to the class $\Tt^{\nu}$. Therefore one has
\begin{lemma}\label{PHsd}
Let $P$ given by \eqref{Pd} with $g$ smooth, real  and having a zero of order at least 3 at the origin. Then $P\in \Hs\cap \Tt^{\nu}$ for all $s>d/2$ and for $\nu>d/8$.
\end{lemma}
We also need a $d$-dimensional definition of {\em normal form} homogeneous polynomial :
\begin{definition}\label{def:NFd}
Let $k = 2m$ be an even integer, a formal polynomial $Z$ homogeneous of degree $k$ on $\Qs$ is in normal form if it reads
$$Z(\xi,\eta)=\sum_{j\in \Nd^k}\quad\sum_{l_{1},l'_1=1}^{d_{j_{1}}}\ldots\sum_{l_{k},l'_k=1}^{d_{j_{k}}}
a_{j,l,l'}\xi_{j_{1},l_{1}}\eta_{j_{1},l'_{1}}\ldots \xi_{j_{k},l_{k}}\eta_{j_{k},l'_{k}}$$
for all $(\xi,\eta)\in \Qs$. \end{definition} 
One easily verifies that if $Z$ is in normal form then $Z$ commutes with each $J_j=\sum_{l=1}^{d_{j}}\xi_{j,l}\eta_{j,l} $ since for instance
$$\{ \xi_{j_{1},l_{1}}\eta_{j_{1},l'_{1}}, \xi_{j_{1},l_{1}}\eta_{j_{1},l_{1}}+\xi_{j_{1},l'_{1}}\eta_{j_{1},l'_{1}}\}=0.$$
Modifying slightly the proof of Theorem \ref{thm:birk} we get
\begin{theorem}\label{birk-gene}
    Let $P$ be a real Hamiltonian belonging in  $\Tt^{\nu}\cap \Hsd$  for some $\nu \geq0$ and for all $s$ sufficiently large and let  $\omega$ be a weakly 
   non resonant frequency vector in the sense of \eqref{A.2d}. Then 
    for any $r\geq 3$ there exists $s_{0}$ and  for any $s\geq 
    s_{0}$ there exists $\U$, $\V$ neighborhoods of the origin in 
    $\Qs$ and 
    $\tau : \V \to \U$ a real analytic canonical transformation  
     which puts $H=H_{0}+P$ in normal form up 
    to order $r$ i.e.
    $$H\circ \tau= H_{0}+Z+R$$
   
    with
       
    \bi
    \item[(i)] $Z$ is a real continuous polynomial of degree $r$ which belongs to $\Hsd$ and which is in normal form in the sense of Definition \ref{def:NFd}. In particular $Z$
    commutes with all  $J_{j}$, $j\geq 1$, i.e. $\{Z, 
    J_{j}\}=0$ for all $j\geq 1$.
     \item[(ii)]  $R$ is real and belongs to $\Hsd$, furthermore $\norma{X_{R}(z)}_{s}\leq 
    C_{s}\norma{z}_{s}^{r}$ for all $z\in \V_{s}$.
    \item[(iii)] $\tau$ is close to the identity: $\norma{\tau 
    (z)-z}_{s}\leq C_{s}\norma{z}_{s}^2$ for all $z\in \V_{s}$.
    \ei
\end{theorem}
\proof The only new point when comparing with Theorem \ref{thm:birk},
is that  in assertion (ii) we obtain $\{Z, 
    J_{j}\}=0$ for all $j\geq 1$ instead of $\{Z, 
    I_{j}\}=0$ for all $j\geq 1$. 
 Actually, in view of \eqref{A.2d}, we adapt  Lemma \ref{homo}, and in particular \eqref{bc1} and \eqref{bc2},  in such a way $\chi \in \Tt^{\nu,+}$  and  $Z$ is in normal form in the sense of Definition \ref{def:NFd}.  
  
On the other hand, we also verify, following the lines of the proof of assertion (iv) of Proposition \ref{tame}, that a  homogeneous polynomial of degree $k+1$ in normal form $Z\in \Tt^\nu$ satisfies $\norma{X_Z(z)}_s\leq C\norma{z}^k_s$ for all $z$ in a neighborhood of the origin. In particular, if $Z\in \Tt^\nu$  is in normal form, it automatically belongs to $\Hsd$ (this point was crucial in the proof of Theorem \ref{thm:birk}).\\
 \endproof

Notice that the normal form $H_{0}+Z$ is no longer, in general, 
integrable. The dynamical consequences are the same as in Theorem 
\ref{thm:dyn} (i) and (ii) but we have to replace $I_{j}$ by $J_{j}$ in the 
second assertion. Actually 
the
$J_{j}$ play the r\^ole of almost actions: they are almost conserved 
quantities.

\begin{theorem}\label{thm:dynd} 
 Assume that  $m\in F_{k}$ defined in  proposition \ref{resd} and that $g$ is $ C^\infty $ on a neighborhood of $0$ in $\C^2$, $g$ is real i.e. $g(z,\bar z)\in \R$ and $g$ vanishes at least at order 3 at the origin. For each $r\geq 3$ and 
 $s\geq s_{0}(r)$, there exists 
 $\e_{0}>0$ and $c>0$ such that  for any  $\psi_0$ in $\tHs$, any 
$\epsilon \in (0,\epsilon_{0})$, the equation
$$i\psi_t=(-\Delta +x^{2} +M)\psi +\partial_2 g(\psi,\bar \psi), \quad x\in \R^d,\  t\in \R$$ 
 with Cauchy data $\psi_0$
has a unique solution 
$ \psi\in 
C^1((-T_{\epsilon},T_{\epsilon}),\tilde H^{s})$ with
$$T_{\epsilon}\geq c\epsilon^{-r}.$$ Moreover 
for any $t\in (-T_{\epsilon},T_{\epsilon})$, one has
$$\Vert{\psi(t,\cdot )}\Vert_{\tilde H^{s}}\leq 2\epsilon$$ and
$$\sum_{j\geq 1}j^{s}|J_j(t)-J_j(0)|\leq {\e^{3}} $$
where $J_j(t)=\sum_{l=1}^{d_{j}}\va{\xi_{j,l}}^2$, $j\geq 1$ are the "pseudo-actions" of $\psi(t,\cdot)=\sum_{j,l} \xi_{j,l}(t)\Phi_{j,l}(\cdot)$ .
\end{theorem}
\proof
Just remark that as in the proof of Theorem \ref{thm:dyn},  defining
$N(z):=2\sum_{j\in \Nd} j^{s}J_j(\xi,\eta)=2\sum_{j\in \Nd} j^{s}\sum_{l=1}^{d_j} \xi_{j,l}\eta_{j,l}$ one has $N(z)=\norma{z}_{s/2}^2$ for all real point $z=(\xi,\bar\xi)$. On the other hand, using that $Z$ commutes with $J_j$, we have
$$
\{N\circ \tau^{-1},H\}(z)=\{N ,H\circ \tau \}\circ \tau^{-1}(z)= \{N, 
R\}(z').$$ 
Therefore,  in the normalized 
variables, we have the estimate
$|\dot N|\leq C N^{(r+1)/2 }$ and the theorem follows as in the proof of Theorem \ref{thm:dyn}.

\endproof
\bibliographystyle{amsalpha}

\bigskip

{\bf Beno\^{i}t Gr\'ebert, Rafik Imekraz, \'Eric Paturel}

{\it Laboratoire de Math\'ematiques Jean Leray UMR 6629,

Universit\'e de Nantes,

2, rue de la Houssini\`ere,

44322 Nantes Cedex 3, France}

\bigskip

E-mail: \quad
\parbox{5cm} {\tt benoit.grebert@univ-nantes.fr\\ rafik.imekraz@univ-nantes.fr\\ eric.paturel@univ-nantes.fr}

\end{document}